\documentclass[a4paper,10pt]{article}
\usepackage{amsmath,amsthm,amssymb,vmargin}
\usepackage{graphicx}
\allowdisplaybreaks
\title{Hitting spheres on hyperbolic spaces}
\date{}

\setmarginsrb{25mm}{20mm}{25mm}{39mm}{0mm}{20mm}{0mm}{10mm}

\newtheorem{teo}{Theorem}[section]

\newtheorem{lem}{Lemma}[section]
\newtheorem{cor}{Corollary}[section]
\theoremstyle{definition}
\newtheorem{rem}{Remark}[section]

\numberwithin{equation}{section}

\newcommand{\Dim}{\noindent\textsc{Proof\\}}
\newcommand{\EndDim}{ $\blacksquare$ \\}

\newcommand{\bm}[1]{\mbox{\boldmath{$#1$}}}

\begin{document}
\maketitle

\long\def\symbolfootnote[#1]#2{\begingroup
\def\thefootnote{\fnsymbol{footnote}}\footnote[#1]{#2}\endgroup} 
\author{
{\center
{ \large Valentina Cammarota \symbolfootnote[2]{Dipartimento di Statistica, Sapienza Universit\`a di Roma,  P.le Aldo Moro 5, 00185 Rome, Italy.  Tel.: +390649910499, fax: +39064959241.  E-mail address: valentina.cammarota@uniroma1.it.}\;\;\;\;\;\;  Enzo Orsingher \symbolfootnote[9]{{\it Corresponding author}. Dipartimento di Statistica, Sapienza Universit\`a di Roma,  P.le Aldo Moro 5, 00185 Rome, Italy.  Tel.: +390649910585, fax: +39064959241. E-mail address: enzo.orsingher@uniroma1.it. }} \\
}
}
\vspace{1cm}

\begin{abstract}
For a hyperbolic Brownian motion on the Poincar\'e half-plane $\mathbb{H}^2$, starting from a point of hyperbolic coordinates $z=(\eta, \alpha)$ inside a hyperbolic disc $U$ of radius $\bar{\eta}$, we obtain the probability of hitting the boundary $\partial U$ at the point $(\bar \eta,\bar \alpha)$. For $\bar{\eta} \to \infty$ we derive the asymptotic Cauchy hitting distribution on $\partial \mathbb{H}^2$ and for small values of $\eta$ and $\bar \eta$ we obtain the classical  Euclidean Poisson kernel. The exit probabilities $\mathbb{P}_z\{T_{\eta_1}<T_{\eta_2}\}$ from a hyperbolic annulus in $\mathbb{H}^2$ of radii $\eta_1$ and $\eta_2$ are derived and the transient behaviour of hyperbolic Brownian motion is considered. Similar probabilities are calculated also for a Brownian motion on the surface of the three dimensional sphere.

For the hyperbolic half-space $\mathbb{H}^n$ we obtain the Poisson kernel of a ball in terms of a series involving Gegenbauer polynomials and hypergeometric functions. For small domains in $\mathbb{H}^n$ we obtain the $n$-dimensional Euclidean Poisson kernel. The exit probabilities from an annulus are derived also in the $n$-dimensional case.
\end{abstract}

\begin{small}

\end{small}


\section{Introduction}

Hyperbolic Brownian motion has been studied over the years by several authors on the half-plane $\mathbb{H}^2$ and on the Poincar\'e disc $\mathbb{D}^2$ and more recently in the $n$-dimensional hyperbolic space (see, for example, Matsumoto and Yor \cite{matsumoto}, Gruet \cite{gruet}, Byczkowski et al. \cite{byczkowski} and Byczkowski and Malecki \cite{byczkowski2}). The hyperbolic half-space $\mathbb{H}^n$ is given by
$$\mathbb{H}^n=\{z=(x,y): x \in \mathbb{R}^{n-1}, y>0\}$$
with the distance formula 
\begin{equation*}
\cosh \eta(z',z)=1+\frac{||z'-z||^2}{2 y y'}.
\end{equation*}
The hyperbolic Brownian motion is a diffusion governed by the generator 
\begin{equation}  \label{subj}
\Delta_n=\frac{y^2}2 \left( \sum_{i=1}^{n-1} \frac{\partial^2}{\partial x_i^2}+\frac{\partial^2}{\partial y^2}\right)-\frac{(n-2)}2 y \frac{\partial}{\partial y}
\end{equation}
(see, for example, Gruet \cite{gruet}). Therefore the probability density $p(x_1,\dots,x_{n-1},y,t)$ of hyperbolic Brownian motion is solution to the Cauchy problem 
\begin{equation*}
\frac{\partial p}{\partial t}=\frac{y^2}2 \left( \sum_{i=1}^{n-1} \frac{\partial^2 p}{\partial x_i^2}+\frac{\partial^2 p}{\partial y^2}\right)-\frac{(n-2)}2 y \frac{\partial p}{\partial y}
\end{equation*}
subject to the initial condition 
\begin{equation*}
p(x_1,\dots,x_{n-1},y,0)=\prod_{j=1}^{n-1} \delta(x_j)\; \delta(y-1).
\end{equation*}
For our purposes it is important to express the generator (\ref{subj}) in hyperbolic coordinates $(\eta,\bm{\alpha})=(\eta,\alpha_1,\dots,\alpha_{n-1})$ as follows 
\begin{align} \label{radial}
\Delta_n=\frac{\partial^2}{\partial \eta^2}+\frac{n-1}{\tanh \eta} \frac{\partial}{\partial \eta} + \frac{1}{\sinh^2 \eta} \Delta_{S_{n-1}}
\end{align}
where $\Delta_{S_{n-1}}$ is the Laplace operator on the $(n-1)$-dimensional unit sphere (see, for example, Helgason \cite{helgason} page 158 or Grigor'yan \cite{grigoryan}). 

The aim of this paper is to study the hitting distribution on a hyperbolic sphere for a hyperbolic Brownian motion  starting from an arbitrary point inside the sphere. Our work is related to the paper by Byczkowski et al. \cite{byczkowski} where the Poisson kernel of half-spaces in $\mathbb{H}^n$, $n>2$, is studied and the paper by Byczkowski and Malecki \cite{byczkowski2} where the Poisson kernel of a ball in the  Poincar\'e disc $\mathbb{D}^n$, $n>2$, is considered.

The first part of our paper concerns the derivation of the Poisson kernel of a hyperbolic disc in $\mathbb{H}^2$ by solving the Dirichlet problem 
\begin{equation} \label{giga}
\begin{cases}
\left[ \frac{\partial^2}{\partial \eta^2}+\frac{1}{\tanh \eta} \frac{\partial}{\partial \eta}+\frac{1}{\sinh^2 \eta} \frac{\partial^2}{\partial \alpha^2}  \right] u (\eta,\alpha; \bar{\eta}, \bar{\alpha})=0, &0<\eta < \bar{\eta}<\infty,\\
u(\bar{\eta}, {\alpha}; \bar{\eta},\bar{\alpha})=\delta(\alpha-\bar{\alpha}), &  {\alpha}, \bar{\alpha} \in (-\pi, \pi].
\end{cases}
\end{equation}
The interplay between Dirichlet problems and hitting probabilities in various contexts is outlined, for example, in Grigor'yan \cite{grigoryan}. The explicit solution of (\ref{giga}) is 
\begin{align} \label{cocca}
u( \eta,\alpha;\bar{\eta}, \bar{\alpha})=\frac{1}{2 \pi}\frac{\cosh \bar{\eta}-\cosh \eta}{\cosh \eta \cosh \bar{\eta}-1-\sinh \eta \sinh \bar{\eta} \cos (\alpha-\bar{\alpha})}
\end{align}
and represents the hitting distribution on the hyperbolic circumference of radius $\bar{\eta}$ for the hyperbolic Brownian motion starting at $(\eta,\alpha)$.

We show that for $\bar{\eta}\to \infty$ the distribution (\ref{cocca}) tends to the Cauchy distribution as was found by means of other arguments in Baldi et al. \cite{casadio}.

The solution to the Dirichlet problem (\ref{giga}) is carried out by two different approaches. One is based on the direct solution of the hyperbolic Laplace equation and the second one is based on some integral representation of the associated Legendre polynomials. 

The derivation of the $n$-dimensional Poisson kernel for $n>2$ is much more clumsy and the final expression is given as a series involving Gegenbauer polynomials and hypergeometric functions. A substantial simplification in the calculations is obtained by applying a suitable rotation of the hyperbolic sphere so that the kernel can be expressed in terms of the hyperbolic distance $\eta$ and the angle $\alpha_1-\bar{\alpha}_1$ between the geodesic lines with ends points $(\eta,\bm \alpha)$ and $(\bar{\eta},{\bar {\bm \alpha}})$. Its explicit form reads 
\begin{align} \label{fif}
&u(\eta, \alpha_1; \bar{\eta}, \bar{\alpha}_1) \nonumber \\
&\;=\frac{\Omega_{n-1}}{\Omega_n}   \sum_{k=0}^\infty \left(\frac{2k}{n-2}+1 \right)   \frac{  \tanh^k \frac{{\eta}} 2\; F\left(k,1-\frac n 2;k+\frac n 2;  \tanh^2 \frac{ {\eta}} 2 \right)}{ \tanh^k \frac{ \bar{\eta}} 2\; F\left(k,1-\frac n 2;k+\frac n 2;  \tanh^2 \frac{ \bar{\eta}} 2 \right)} C^{(\frac{n-2}{2})}_k(\cos( \alpha_1-\bar{\alpha}_1)) \sin^{n-2} (\alpha_1-\bar{\alpha}_1),
\end{align}
where $n>2$, $0<\eta < \bar{\eta}<\infty$, $\alpha_1-\bar{\alpha}_1 \in (0, \pi]$, $\Omega_n=\frac{2 \pi^{\frac{n}{2}}}{\Gamma(\frac{n}{2})}$ is the surface area of the $n$-dimensional Euclidean unit sphere, $F(\alpha,\beta;\gamma,x)$ is the hypergeometric function and $C_k^{(n)}(x)$ are the Gegenbauer polynomials.  

The Poisson kernel of a ball in the hyperbolic disc $\mathbb{D}^n$, $n>2$, is obtained in Byczkowski and Malecki \cite{byczkowski2}, formula (16) and must be compared with (\ref{fif}) above.

Unfortunately formula (\ref{fif}) cannot be reduced to a fine form as (\ref{cocca}). However, for sufficiently small domains, we extract from (\ref{fif}) the $n$-dimensional  Euclidean  Poisson kernel.

Section 3 is devoted to the exit probabilities  $\mathbb{P}_z\{T_{\eta_1}<T_{\eta_2}\}$ from a hyperbolic annulus of radii $\eta_1$ and $\eta_2$. We examine in detail both the planar and the higher dimensional case discussing also the transient behaviour of hyperbolic Brownian motion. 

In the last section the hitting probabilities on a spherical circle for a spherical Brownian motion starting from $p=(\vartheta,\varphi)$ are considered. In particular the most interesting result here is that 
$$\mathbb{P}_p\{B_{S}(T_{\bar{\vartheta}}) \in \mathrm{d} {\bar{\varphi}} \}=\frac{1}{2 \pi} \frac{\cos \vartheta-\cos {\bar \vartheta} }{1-\cos \theta \cos \bar{\vartheta} - \sin \vartheta \sin{\bar \vartheta} \cos(\varphi-\bar{\varphi})}\mathrm{d} \bar {\varphi},\hspace{0.5cm}0<\bar \vartheta<\vartheta<\pi,\;\; {\varphi}, \bar{\varphi} \in (0, 2 \pi].$$

\section{Hitting distribution on a hyperbolic sphere in $\mathbb{H}^n$}

\subsection{Two dimensional case}

We study here the Poisson kernel of the circle in the hyperbolic plane $\mathbb{H}^2=\{(x,y):x \in \mathbb{R}, y>0\}$ endowed with the Riemannian metric 
$$\mathrm{d}s^2=\frac{\mathrm{d}x^2+\mathrm{d}y^2}{y^2},$$
and the distance formula 
\begin{equation} \label{ciccio}
\cosh \eta(z',z)=\frac{(x'-x)^2+y'^2+y^2}{2 y y'}.
\end{equation}
We denote with $\eta$ the hyperbolic distance from the origin $O=(0,1)$ of $\mathbb{H}^2$. The Laplace operator on $\mathbb{H}^2$ in cartesian coordinates reads 
\begin{equation} \label{1.5}
\Delta=y^2 \left(\frac{\partial^2}{\partial x^2}+\frac{\partial^2}{\partial y^2}\right)
\end{equation}
(for a proof see, for example, Chavel \cite{chavel} page 265). It is convenient to write the Laplace operator  in hyperbolic coordinates $(\eta, \alpha)$ 
\begin{equation} \label{lao}
\Delta= \frac{\partial^2}{\partial \eta^2}+\frac{1}{\tanh \eta} \frac{\partial}{\partial \eta}+\frac{1}{\sinh^2 \eta} \frac{\partial^2}{\partial \alpha^2}
\end{equation}
(for information on hyperbolic coordinates see Cammarota and Orsingher \cite{travel}). The relationship between hyperbolic coordinates $(\eta, \alpha)$ and the cartesian coordinates $(x,y)$ is given by
\begin{equation} \label{hjhj}
\begin{cases}
{x=\frac{\sinh \eta \cos \alpha}{\cosh \eta-\sinh \eta \sin \alpha}}, \\
{y=\frac{1}{\cosh \eta-\sinh \eta \sin \alpha}}.
\end{cases}
\end{equation}    
By exploiting (\ref{hjhj}), in the paper by  Lao and Orsingher \cite{lao}, the Laplace operator (\ref{lao}) is obtained from (\ref{1.5}). \\

\noindent We have now our first theorem.

\begin{teo} \label{theteo}
Let $U=\{(\eta,\alpha): \eta<\bar{\eta}\}$ be a hyperbolic disc in $\mathbb{H}^2$ with radius $\bar{\eta}$ and center in $O$, the solution to the Dirichlet problem 
\begin{equation} \label{dsds}
\begin{cases}
\left[ \frac{\partial^2}{\partial \eta^2}+\frac{1}{\tanh \eta} \frac{\partial}{\partial \eta}+\frac{1}{\sinh^2 \eta} \frac{\partial^2}{\partial \alpha^2}  \right] u(\eta,\alpha;\bar{\eta}, \bar{\alpha})=0, &0<\eta < \bar{\eta}<\infty,\\
u(\bar{\eta},\alpha;\bar{\eta}, \bar{\alpha})=\delta(\alpha-\bar \alpha), &  {\alpha}, \bar{\alpha} \in (-\pi, \pi],
\end{cases}
\end{equation}
is given by 
\begin{eqnarray} \label{statement}
u(\eta,\alpha;\bar{\eta}, \bar{\alpha})=\frac{1}{2 \pi} \frac{\cosh \bar{\eta}-\cosh \eta}{\cosh \eta \cosh \bar{\eta}-1-\sinh \eta \sinh \bar{\eta} \cos (\alpha-\bar{\alpha})}.
\end{eqnarray}
\end{teo}
\Dim 
Our proof is based on the classical method of separation of variables. We assume that 
\begin{equation} \label{28}
u(\eta,\alpha;\bar{\eta}, \bar{\alpha})=E(\eta) \Theta(\alpha)
\end{equation}
and we arrive at the following ordinary equations 
\begin{equation} \label{kkk}
\begin{cases}
\Theta''(\alpha)+\mu^2 \: \Theta(\alpha)=0,\\
\sinh^2 \eta\;  E''(\eta)+\cosh \eta\; \sinh \eta \;E'(\eta)-\mu^2 E(\eta)=0, 
\end{cases}
\end{equation} 
where $\mu^2$ is an arbitrary constant. The first equation has general solution
\begin{equation} \label{s1}
\Theta(\alpha)=A \cos(\mu \alpha)+B \sin(\mu \alpha)
\end{equation}
and becomes periodic with period $2 \pi$ for $\mu=m \in \mathbb{N}$. The second equation necessitates some further treatment. We start with the change of variable $w=\cosh \eta$ which transforms the second equation of (\ref{kkk}) into 
\begin{eqnarray} \label{ss}
(1-w^2) G''(w)-2 w G'(w)-\frac{m^2}{1-w^2}G(w)&=&0. 
\end{eqnarray}
The general solution to (\ref{ss}) can be conveniently written as  
\begin{eqnarray} \label{hhh}
G(w)=C_1 \left| \frac{w+1}{w-1}\right|^{m/2}+C_2 \left| \frac{w-1}{w+1}\right|^{m/2}, \hspace{2cm} m \ne 0,
\end{eqnarray}
(see, for example, Polyanin and Zaitsev  \cite{polyanin} Section 2.1.2, formula 233 for $a=1$, $b=-1$, $\lambda=0$ and $\mu=-m^2$). From (\ref{hhh}) we have that 
\begin{align}\label{ggg}
E(\eta)=C_1 \left( \frac{\cosh \eta+1}{\cosh \eta-1}\right)^{m/2}+C_2 \left( \frac{\cosh \eta-1}{\cosh \eta+1}\right)^{m/2}=C_1 \left(\frac{\cosh \eta+1}{\sinh \eta}\right)^m+C_2 \left( \frac{\cosh \eta-1}{\sinh \eta}\right)^m.
\end{align}
We disregard the first term of (\ref{ggg}) since our aim is to extract finite-valued and increasing solutions to (\ref{dsds}), so that we have 
\begin{eqnarray} \label{s2}
E(\eta)=C \left( \frac{\cosh \eta-1}{\sinh \eta}\right)^m=C \tanh^m \frac \eta 2.
\end{eqnarray}
In light of (\ref{28}), (\ref{s1}) and (\ref{s2}) we can write 
\begin{eqnarray} \label{2.9}
u(\eta,\alpha; \bar{\eta}, \bar{\alpha})&=&\sum_{m=0}^\infty  \Theta_m(\alpha)  E_m(\eta) =A_0+\sum_{m=1}^\infty [A_m \cos(m \alpha)+B_m \sin(m \alpha)] \left(\frac{\cosh \eta-1}{\sinh \eta} \right)^m.\;\;
\end{eqnarray}
If we take the Fourier expansion of the Dirac delta function
\begin{eqnarray} \label{2.15}
\delta(\alpha-\bar{\alpha})&=&\frac{1}{2 \pi}+\frac{1}{\pi} \sum_{m=1}^\infty \cos [m (\alpha-\bar{\alpha})]\nonumber \\
&=&\frac{1}{2 \pi}+\frac{1}{\pi} \sum_{m=1}^\infty [ \cos (m \alpha) \cos (m \bar{\alpha}) + \sin (m \alpha) \sin (m \bar{\alpha}) ],
\end{eqnarray}
by comparing (\ref{2.9}) with (\ref{2.15}) we obtain the Fourier coefficients $A_m$ and $B_m$ so that we can write
\begin{eqnarray} \label{ea}
u(\eta,\alpha;\bar{\eta}, \bar{\alpha})&=&\frac{1}{2 \pi}+\frac{1}{\pi} \sum_{m=1}^\infty [ \cos (m \alpha) \cos (m \bar{\alpha})+\sin (m \alpha) \sin (m \bar{\alpha}) ]\left(\frac{\cosh \bar{\eta}-1}{\sinh \bar{\eta}} \right)^{-m}  \left(\frac{\cosh \eta-1}{\sinh \eta} \right)^m \nonumber\\
&=&\frac{1}{2 \pi}+\frac{1}{\pi} \sum_{m=1}^\infty  \cos (m (\alpha-\bar{\alpha})) \left(\frac{\cosh \bar{\eta}-1}{\sinh \bar{\eta}} \right)^{-m}  \left(\frac{\cosh \eta-1}{\sinh \eta} \right)^m \nonumber \\
&=&\frac{1}{2\pi} \left[ 1+   \sum_{m=1}^\infty  \left[ \left(e^{i (\alpha-\bar{\alpha})} \frac{\sinh \bar{\eta}}{\cosh \bar{\eta}-1} \frac{\cosh \eta-1}{\sinh \eta}     \right)^m + \left(e^{-i (\alpha-\bar{\alpha})} \frac{\sinh \bar{\eta}}{\cosh \bar{\eta}-1}  \frac{\cosh \eta-1}{\sinh \eta}    \right)^m\right] \right] \nonumber\\
&=& \frac{1}{2 \pi}  \frac{\left(\frac{\cosh \bar{\eta}-1}{\sinh \bar{\eta}}   \right)^2-  \left(\frac{\cosh \eta-1}{\sinh \eta}   \right)^2}{ \left(\frac{\cosh \bar{\eta}-1}{\sinh \bar{\eta}}   \right)^2+ \left(\frac{\cosh \eta-1}{\sinh \eta}   \right)^2-2  \frac{\cosh\bar{ \eta}-1}{\sinh \bar{\eta}}   \frac{\cosh \eta-1}{\sinh \eta}    \cos (\alpha-\bar{\alpha})}\\
&=&\frac{1}{2 \pi}  \frac{\tanh^2 \frac{\bar \eta}{2}-  \tanh^2 \frac{\eta}{2}}{ \tanh^2 \frac{\bar \eta}{2}+ \tanh^2 \frac{ \eta}{2}-2  \tanh \frac{\bar \eta}{2}  \tanh \frac{ \eta}{2}  \cos (\alpha-\bar{\alpha})} \nonumber. 
\end{eqnarray}
The expression in (\ref{ea}) can be substantially simplified by observing that:  
\begin{align*}
&(\cosh \bar{\eta}-1)^2 \sinh^2 \eta-(\cosh \eta-1)^2 \sinh^2 \bar{\eta}=(\cosh \bar{\eta}-1) (\cosh {\eta}-1)[2 \cosh \bar{\eta}-2 \cosh \eta]
\end{align*}
and
\begin{align*}
&(\cosh \bar{\eta}-1)^2 \sinh^2 \eta+ (\cosh {\eta}-1)^2 \sinh^2 \bar{\eta}-2 (\cosh \bar{\eta}-1) (\cosh {\eta}-1) \sinh \bar{\eta} \sinh \eta \cos (\alpha-\bar{\alpha})\\
&\hspace{2cm}=(\cosh \bar{\eta}-1) (\cosh {\eta}-1)[2 \cosh \eta \cosh \bar{\eta}  -2-2 \sinh \eta \sinh \bar{\eta} \cos (\alpha-\bar{\alpha})].
\end{align*}
In view of all these calculations we have that the hyperbolic Poisson kernel takes the form 
\begin{equation*}
u( \eta,\alpha;\bar{\eta}, \bar{\alpha})=\frac{1}{2 \pi}\frac{\cosh \bar{\eta}-\cosh \eta}{\cosh \eta \cosh \bar{\eta}-1-\sinh \eta \sinh \bar{\eta} \cos (\alpha-\bar{\alpha})}.
\end{equation*} 
\EndDim

\begin{rem} \label{therem}
It is possible to obtain the expression (\ref{2.9}) by means of an alternative approach as follows. We start from the associated Legendre equation  
\begin{eqnarray} \label{leg}
(1-z^2) y''(z)-2zy'(z)+ \left[ \nu(\nu+1)-\frac{m^2}{1-z^2} \right]y(z)=0
\end{eqnarray}
which coincides with (\ref{ss}) for $\nu=0$ or $\nu=-1$. In view of Gradshteyn and Ryzhik \cite{G-R}  formula 8.711.2, the solution to (\ref{leg}) can be written as 
\begin{equation*}
P_{\nu}^m (z)= \frac{(-1)^m}{\pi} \frac{\Gamma(\nu+ 1)}{\Gamma(\nu-m+1)}  \int_0^{\pi}  \frac{\cos (m \varphi)}{(z+\sqrt{z^2-1} \cos \varphi)^{\nu+ 1}} \mathrm{d}\varphi, \hspace{1cm} |\arg z|<\frac \pi 2.
\end{equation*}
If $\nu=-1$ we have
\begin{eqnarray*}
P_{-1}^m(\cosh \eta)= \frac{(-1)^m}{ \pi} \frac{1}{\Gamma(-m)}  \int_0^{ \pi} \cos (m \phi) \mathrm{d}\phi=0, \hspace{1cm} \mathrm{for}\;\; m \in \mathbb{Z}.
\end{eqnarray*}
If $\nu=0$ we have 
\begin{eqnarray*} 
P_{0}^m(\cosh \eta)= \frac{(-1)^m}{ \pi} \frac{1}{\Gamma(1-m)}  \int_0^{ \pi}  \frac{\cos (m \phi)}{\cosh \eta+\sinh \eta \cos \phi} \mathrm{d}\phi
\begin{cases}
=0,  \hspace{1cm} \mathrm{for}\;\; m=1,2,\dots \\
\ne 0,  \hspace{1cm} \mathrm{for}\;\; m=0,-1,-2,\dots
\end{cases} 
\end{eqnarray*}
It follows that 
\begin{eqnarray} \label{onde}
u(\eta,\alpha;\bar{\eta}, \bar{\alpha})&=&\sum_{m=-\infty}^0  [A_m \cos(m \alpha)+B_m \sin(m \alpha)] P_0^m(\cosh \eta) \nonumber \\
&=&\frac{1}{ \pi} \sum_{m=-\infty}^0  \frac{(-1)^m}{ \Gamma(1-m)}  [A_m \cos(m \alpha)+B_m \sin(m \alpha)]  \int_0^{ \pi}  \frac{\cos (m \phi)}{\cosh \eta+\sinh \eta \cos \phi} \mathrm{d}\phi \nonumber \\
&=& \frac{1}{ \pi}\sum_{m=0}^\infty (-1)^m [A_m \cos(m \alpha)+B_m \sin(m \alpha)]  \int_0^{ \pi}  \frac{\cos (m \phi)}{\cosh \eta+\sinh \eta \cos \phi} \mathrm{d}\phi,
\end{eqnarray}
in the last step $A_n$ and $B_n$ include the multiplicative constant $\frac{1}{ \Gamma(1-m)} $. Since we have
\begin{equation}  \label{pp}
1+ 2  \sum_{n=1}^\infty \left(\frac{1-\cosh \eta}{\sinh \eta}   \right)^n  \cos (n \phi)=\frac{1}{\cosh \eta+ \sinh \eta \cos \phi},
\end{equation}
by inserting (\ref{pp}) into (\ref{onde}) we get 
\begin{eqnarray*}
&&\hspace{-1.0cm}u(\eta,\alpha;\bar{\eta}, \bar{\alpha})\\
&=& \frac{1}{ \pi}\sum_{m=0}^\infty (-1)^m [A_m \cos(m \alpha)+B_m \sin(m \alpha)]  \int_0^{ \pi} \cos (m \phi) \left(  1+ 2  \sum_{n=1}^\infty \left(\frac{1-\cosh \eta}{\sinh \eta}   \right)^n  \cos (n \phi) \right)    \mathrm{d}\phi\\
&=&\frac{2}{\pi} \sum_{m=0}^\infty (-1)^m [A_m \cos(m \alpha)+B_m \sin(m \alpha)]    \sum_{n=1}^\infty  \left(\frac{1-\cosh \eta}{\sinh \eta}   \right)^n \int_0^{ \pi} \cos (m \phi)    \cos (n \phi)    \mathrm{d}\phi\\
&=&A_0+\sum_{m=1}^\infty  [A_m \cos(m \alpha)+B_m \sin(m \alpha)]    \left(\frac{\cosh \eta-1}{\sinh \eta}   \right)^m
\end{eqnarray*}
and thus we retrieve (\ref{2.9}).
\end{rem}

\begin{figure}[h] 
 \centering
     {\includegraphics[width=6.5cm, height=6.5cm]{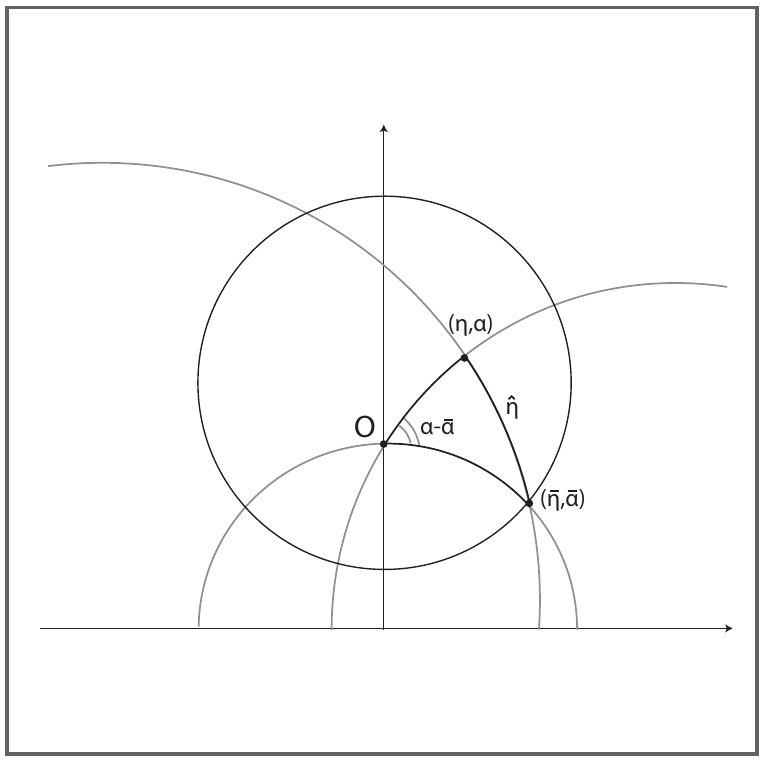}}
 \caption{Hyperbolic triangle in $\mathbb{H}^2$ with sides of length $\eta$, $\bar{\eta}$ and $\hat{\eta}$. } \label{carnottt}
 \end{figure}

\begin{rem}
By applying the hyperbolic Carnot formula we note that it is possible to write the hyperbolic Poisson kernel (\ref{statement}) in a new form. We construct a hyperbolic triangle with sides of length $\eta$, $\bar{\eta}$ and $\hat{\eta}$, and angle between the two sides of length $\eta$ and $\bar{\eta}$ equal to $\theta=\alpha-\bar{\alpha}$, see Figure \ref{carnottt}. The hyperbolic Carnot formula  
$$\cosh \hat{\eta}=\cosh \eta \cosh \bar{\eta}-\sinh \eta \sinh \bar{\eta} \cos (\alpha-\bar{\alpha}),$$
permits us to write (\ref{statement}) as
\begin{equation} \label{uhhu}
u(\eta,\alpha;\bar{\eta}, \bar{\alpha})=\frac{1}{2 \pi} \frac{\cosh \bar{\eta}-\cosh \eta}{\cosh \hat{\eta}-1},
\end{equation}
where the dependence of $u$ from $\alpha$ and $\bar{\alpha}$ is hidden in $\hat{\eta}$.
\end{rem}

\begin{rem}
We observe that the hyperbolic Poisson kernel (\ref{statement}) is a proper probability law. In fact: 
\begin{itemize}
\item It is non-negative because, for $\eta > \bar{\eta}$, we have $\cosh \bar{\eta}-\cosh \eta>0$ and by the hyperbolic Carnot formula 
$$\cosh \eta \cosh \bar{\eta}-1-\sinh \eta \sinh \bar{\eta} \cos (\alpha-\bar{\alpha})=\cosh \hat{\eta}-1>0.$$
\item It integrates to one since it is well-known that \begin{equation}\label{omm} \int_0^{2\pi}\frac{\mathrm{d}\theta}{a+b \cos \theta}=\frac{2 \pi}{\sqrt{a^2-b^2}}\end{equation}
where, in this case, $a=\cosh \eta \cosh \bar{\eta}-1$ and $b=-\sinh \eta \sinh \bar{\eta}$.
\end{itemize}
\end{rem}

\begin{figure}[h]
\begin{minipage}[b]{7.85cm}
\centering
    \includegraphics[width=6.5cm, height=6.5cm]{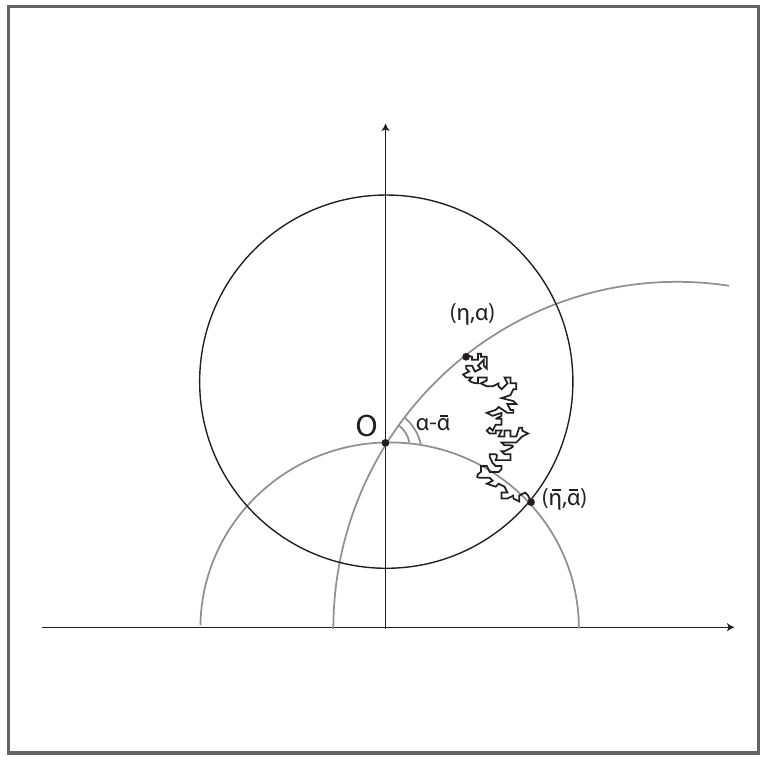}
    \caption{Brownian motion on $\mathbb{H}^2$ starting at $(\eta,\alpha)$ and hitting the boundary of the hyperbolic disc $U$.} \label{treno}
\end{minipage}
\hspace{0.25cm}
\begin{minipage}[b]{7.85cm}
\centering
    \includegraphics[width=6.5cm, height=6.5cm]{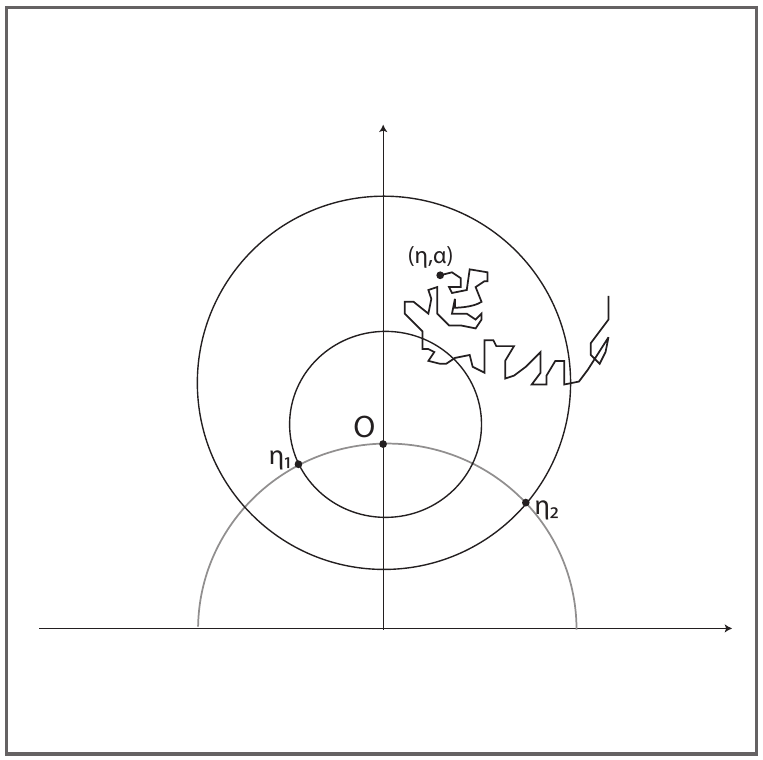}
    \caption{Hyperbolic Brownian motion starting inside the hyperbolic annulus $A$ with radii $\eta_1$ and $\eta_2$.} \label{autobus}
\end{minipage}
 \end{figure}

\begin{rem}
The kernel appearing in formulas (\ref{statement}) and (\ref{uhhu}) represents the law of the position occupied by the hyperbolic Brownian motion $\{B_{\mathbb{H}^2}(t): t\ge0\}$ on $\mathbb{H}^2$ starting from $z=(\eta,\alpha)\in \mathbb{H}^2$ when it hits for the first time the boundary $\partial U$ of the hyperbolic disc $U$. In other words 
\begin{eqnarray*} \label{density}
\mathbb{P}_z\{B_{\mathbb{H}^2}(T_{\bar{\eta}})\in \mathrm{d} \bar{\alpha}\}&=&\frac{1}{2 \pi}\frac{\cosh \bar{\eta}-\cosh \eta}{\cosh \eta \cosh \bar{\eta}-1-\sinh \eta \sinh \bar{\eta} \cos (\alpha-\bar{\alpha})} \mathrm{d}\bar{\alpha},\hspace{0.7cm} \bar{\alpha} \in [0, 2 \pi),
\end{eqnarray*}
where $T_{\bar{\eta}}=\inf\{t>0: B_{\mathbb{H}^2}(t)\in  \partial U \}$, see Figure \ref{treno}.
\end{rem}

\begin{figure}[h]
\begin{minipage}[b]{7.85cm}
\centering
    \includegraphics[width=6.5cm, height=6.5cm]{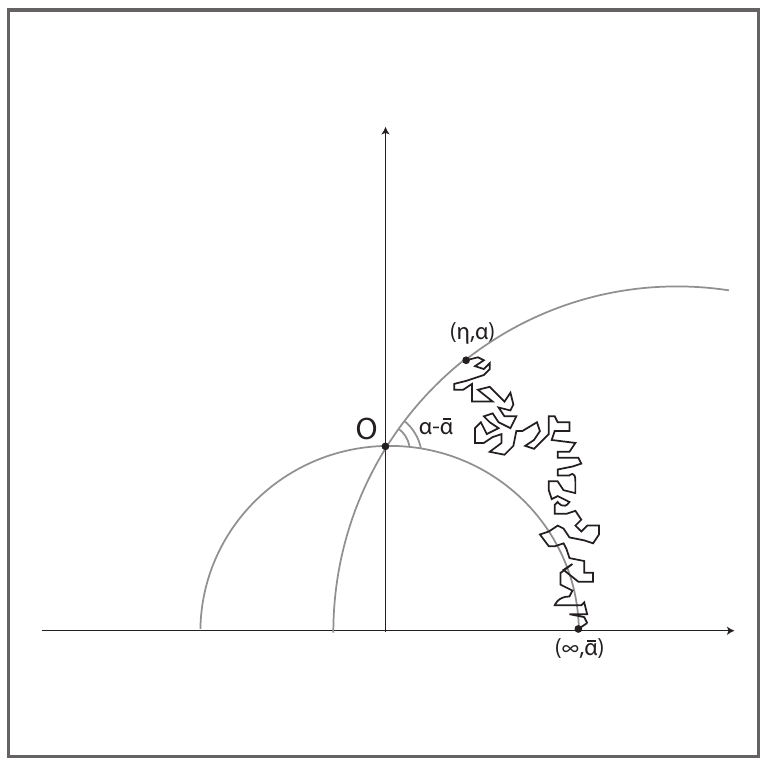} 
    \caption{Brownian motion on $\mathbb{H}^2$ starting at $(\eta,\alpha)$ and hitting the boundary of the hyperbolic plane.} \label{traballa}
\end{minipage}
\hspace{0.25cm}
\begin{minipage}[b]{7.85cm}
\centering
    \includegraphics[width=6.5cm, height=6.5cm]{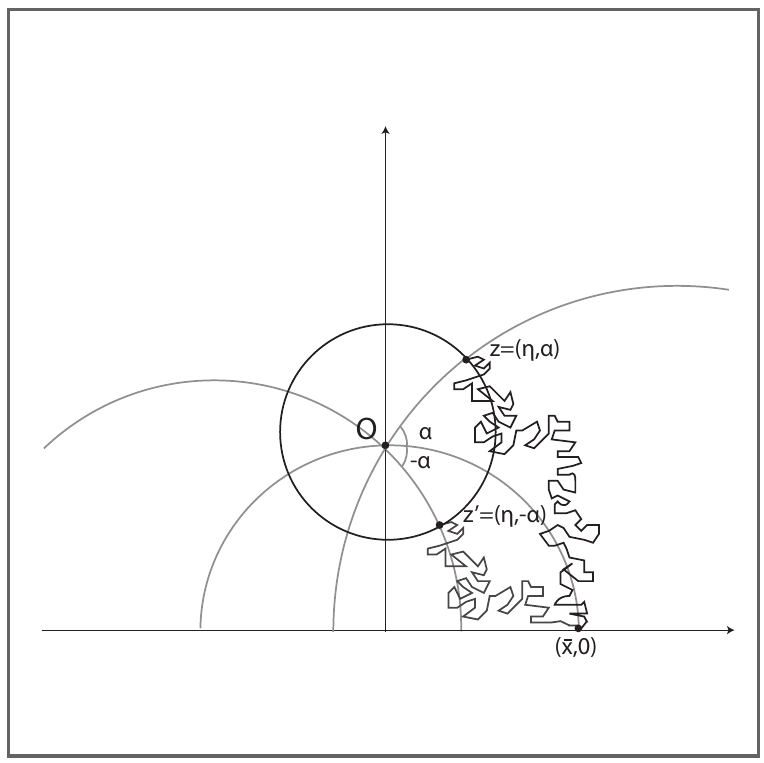}
    \caption{Hyperbolic Brownian motion starting at $z$ and Euclidean Brownian motion starting at $z'$.} \label{tutto}
\end{minipage}
 \end{figure}

\begin{rem}
For small values of $\eta$ and $\bar{\eta}$ the hyperbolic Poisson kernel (\ref{statement}) is approximated by the Euclidean Poisson kernel 
\begin{eqnarray*}
u(\eta,\alpha;\bar{\eta}, \bar{\alpha}) &\sim&\frac{1}{2 \pi} \frac{1+\frac{\bar{\eta}^2}{2} - (1+\frac{{\eta}^2}{2} )}{(1+\frac{\bar{\eta}^2}{2} )(1+\frac{{\eta}^2}{2} )-1- \eta\; \bar{\eta} \cos(\alpha-\bar{\alpha})}= \frac{1}{2 \pi}\frac{\bar{\eta}^2-{\eta}^2}{\bar{\eta}^2+{\eta}^2-2\; {\eta}\; \bar{\eta} \cos(\alpha-\bar{\alpha})}
\end{eqnarray*}
that represents the law of the position occupied by the Euclidean Brownian motion $\{B(t), t \ge 0\}$ on $\mathbb{R}^2$ starting from a point $z=(\eta,\alpha)$ when it hits for the first time the boundary $\partial U$ of the Euclidean disc $U=\{(\eta,\alpha), \eta<\bar{\eta}\}$ with Euclidean radius $\bar{\eta}$. This is a consequence of the fact that in sufficiently small domains of the Lobatchevskian space, the Euclidean geometry is in force. 

\end{rem}
\begin{rem}
We also note that: 
\begin{itemize}
\item For $\eta=0$ formula (\ref{statement}) becomes the uniform distribution as expected.
\item For $\bar{\eta} \to \infty$ we have that 
\begin{equation} \label{lim}
\tilde{u}(\eta,\alpha;\bar{\alpha}):=\lim_{\bar{\eta} \to \infty} u(\eta,\alpha;\bar{\eta}, \bar{\alpha})=\frac{1}{2 \pi} \frac{1}{\cosh \eta-\sinh \eta \cos (\alpha-\bar{\alpha})}.
\end{equation}
\end{itemize}
In view of (\ref{pp}), the limiting distribution (\ref{lim}) can also be written as 
$$\tilde{u}(\eta,\alpha;\bar{\alpha})=\frac{1}{2 \pi} \frac{1}{\cosh \eta-\sinh \eta \cos (\alpha-\bar{\alpha})}=\frac{1}{2 \pi} \left[1+2 \sum_{n=1}^\infty \left(\frac{\cosh \eta-1}{\sinh \eta} \right)^n \cos n(\alpha-\bar{\alpha})\right].$$
We note that $\tilde{u}$ represents the hitting distribution of the hyperbolic Brownian motion, starting at $z=({\eta},{\alpha})$, on the horizontal axis $\partial \mathbb{H}^2=\{(\bar{\eta},\bar{\alpha}): \bar{\eta}=\infty\}=\{(\bar{x},\bar{y}): \bar{y}=0\}$, see Figure \ref{traballa}. We observe that  
the boundary $\partial \mathbb{H}^2$ represents the point at infinity of $\mathbb{H}^2$. We can write the `hitting' probability on $\partial \mathbb{H}^2 $ in the following form
\begin{equation} \label{2.225}
\mathbb{P}_z\{B_{\mathbb{H}^2}(T_{\infty})\in \mathrm{d} \bar{\alpha}\}=\frac{1}{2 \pi} \frac{1}{\cosh \eta-\sinh \eta \cos \alpha \cos \bar{\alpha}-\sinh \eta \sin \alpha \sin \bar{\alpha}} \mathrm{d} \bar{\alpha}.
\end{equation}
We write now  the distribution (\ref{2.225}) in cartesian coordinates. In view of (\ref{hjhj}) we have that 
\begin{equation} \label{uno}
\frac x y=\sinh \eta \cos \alpha,  \hspace{2cm} \tan \alpha=\frac{x^2+y^2-1}{2x}.
\end{equation}
The first relation is an immediate consequence of (\ref{hjhj}) and for a proof of the second equality, see Cammarota and Orsingher \cite{cascade}. From (\ref{ciccio}) and (\ref{uno}) it follows that 
\begin{equation} \label{due}
\sinh \eta \sin \alpha=\sqrt{\cosh^2\eta -1} \frac{\tan \alpha}{\sqrt{1+\tan^2 \alpha}}=\frac{x^2+y^2-1}{2y}.
\end{equation}
Letting ${\bar \eta} \to \infty$ we note, in view of (\ref{hjhj}), that for a point $(\bar{x},\bar{y}) \in \partial \mathbb{H}^2$ it holds that  
\begin{equation} \label{gfd}
\begin{cases}
\bar{x}=\frac{\cos \bar{\alpha}}{1-\sin \bar{\alpha}}, \\
\bar{y}=0.
\end{cases}
\end{equation}    
Formula (\ref{gfd}) implies that $\bar{x}-\cos \bar{\alpha}=\bar{x} \sqrt{1-\cos^2 \bar{\alpha}}$ and this leads to the following relations 
\begin{equation} \label{tre}
\cos \bar{\alpha}=\frac{2 \bar{x}}{1+\bar{x}^2}, \hspace{2cm} \sin \bar{\alpha}=\frac{1-\bar{x}^2}{1+\bar{x}^2}.
\end{equation}
In view of (\ref{uno}), (\ref{due}) and (\ref{tre}) and since $\mathrm{d} \bar{\alpha}=\frac{2}{1+\bar{x}^2} \mathrm{d}\bar{x}$, we can write $\tilde{u}(\eta,\alpha;\bar{\alpha}) \mathrm{d} \bar{\alpha}$ in cartesian coordinates as follows 
\begin{eqnarray} \label{fcv}
\tilde{u}(x,y;\bar{x}) \mathrm{d} \bar{x}&=&\frac{1}{2 \pi} \frac{1}{ \frac{x^2+y^2+1}{2 y}-\frac x y  \frac{2 \bar{x}}{1+\bar{x}^2}-\frac{x^2+y^2-1}{2y} \frac{1-\bar{x}^2}{1+\bar{x}^2} } \frac{2}{1+\bar{x}^2} \mathrm{d} \bar{x}\nonumber \\
&=& \frac{1}{\pi} \frac{2 y}{(x^2+y^2+1) (1+\bar{x}^2)-4x \bar{x}-(x^2+y^2-1)(1-\bar{x}^2)}\mathrm{d} \bar{x}\nonumber \\
&=&\frac{1}{\pi} \frac{y}{\bar{x}^2(x^2+y^2)-2 x \bar{x}+1} \mathrm{d} \bar{x}\nonumber \\
&=& \frac{1}{\pi} \frac{y}{ \left[ \bar{x} \sqrt{x^2+y^2} -\frac{x}{\sqrt{x^2+y^2}}  \right]^2-\frac{x^2}{x^2+y^2}+1} \mathrm{d} \bar{x}\nonumber \\
&=&\frac{1}{\pi} \frac{\frac{y}{x^2+y^2}}{ \left[ \bar{x} -\frac{x}{x^2+y^2}  \right]^2+\left[\frac{y}{x^2+y^2} \right]^2} \mathrm{d} \bar{x}.
\end{eqnarray}

Formula (\ref{fcv}) says that the probability that the hyperbolic Brownian motion starting at $(x,y) \in \mathbb{H}^2$ hits the boundary  of $\mathbb{H}^2$ at $(\bar{x},0)$ is Cauchy distributed with scale parameter $y'=\frac{y}{x^2+y^2}$ and position parameter $x'=\frac{x}{x^2+y^2}$ depending on the starting point. In particular, if the hyperbolic Brownian motion starts at the origin $O$ of $\mathbb{H}^2$, we obtain a standard Cauchy. We note that (\ref{lim}) can be viewed as a Cauchy density in hyperbolic coordinates. 
\end{rem}

\begin{rem}
In view of formula (\ref{fcv}), we also note that the probability that the hyperbolic Brownian motion starting at $z=(x,y)=(\eta,\alpha) \in \mathbb{H}^2$ hits $\partial \mathbb{H}^2$ at $(\bar{x},0)$ is equal to the probability that a Euclidean Brownian motion starting at $z'=(x',y')=(\eta,\alpha')$ hits the $x$-axis at $(\bar{x},0)$, where $z$ and $z'$ have the same hyperbolic distance ${\eta}$ from the origin but  $\alpha'=-\alpha$, see Figure \ref{tutto}. In fact
\begin{align*}
\cosh \eta'&=\frac{\frac{x^2}{(x^2+y^2)^2}+\frac{y^2}{(x^2+y^2)^2}+1}{\frac{2y}{x^2+y^2}}=\frac{x^2+y^2+1}{2y}=\cosh \eta, \\
\tan \alpha'&=\frac{\frac{x^2}{(x^2+y^2)^2}+\frac{y^2}{(x^2+y^2)^2}-1}{\frac{2x}{x^2+y^2}}=\frac{1-x^2-y^2}{2x}=-\tan \alpha.
\end{align*}
Formula (\ref{fcv}) is in accordance with formula (1.2) in Baldi et al. \cite{casadio}. In this paper the hitting distribution  on the horizontal axis, for the hyperbolic Brownian with horizontal and vertical drift,  is obtained from the hitting distribution on the horizontal lines $H_a=\{(x,y)\in \mathbb{H}^2: y=a>0\}$ when $a \to 0$.
\end{rem}

\begin{rem}
The Poisson kernel (\ref{statement}) can be conveniently written also in cartesian coordinates by exploiting the relations (\ref{uno}), (\ref{due}) and the hyperbolic distance formula $$\cosh \eta=\frac{x^2+y^2+1}{2y}.$$ We have that
\begin{align*}
u(x,y;\bar{x}, \bar{y})&=\frac{1}{2 \pi} \frac{\frac{\bar{x}^2+\bar{y}^2+1}{2\bar{y}}-\frac{x^2+y^2+1}{2y}}{\frac{x^2+y^2+1}{2y} \frac{\bar{x}^2+\bar{y}^2+1}{2\bar{y}}-1- \frac x y \frac{\bar{x}}{\bar{y}} -\frac{x^2+y^2-1}{2y} \frac{\bar{x}^2+\bar{y}^2-1}{2\bar{y}}}\\
&=\frac{1}{\pi} \frac{(\bar{x}^2+\bar{y}^2+1)y-({x}^2+{y}^2+1)\bar y}{x^2+y^2+\bar x^2+\bar y^2-2 y \bar y-2 x \bar x} \\
&=\frac{1}{\pi} \frac{(\bar{x}^2+\bar{y}^2)y-({x}^2+{y}^2)\bar y +y-\bar y}{(x-\bar x)^2+(y-\bar y)^2}. 
\end{align*}
In the special case where $\bar y=0$ the previous expression becomes 
\begin{align*}
u(x,y;\bar{x}, 0)=\frac{1}{\pi} \frac{(1+\bar x^2)y}{(x-\bar x)^2+ y^2}
\end{align*}
and thus multiplying by $\frac{1}{1+\bar x^2}$ we get the Cauchy density as expected.
\end{rem}

\subsection{Multidimensional case}

Let $\mathbb{H}^n=\{z=(x,y): x \in \mathbb{R}^{n-1}, y>0\}$ be the $n$-dimensional hyperbolic plane, $n>2$,  with origin $O=(0,\dots,0,1)$ endowed with the Riemannian metric 
\begin{equation*} \label{rim}
\mathrm{d}s^2=\frac{\mathrm{d} x_1^2+\cdots+\mathrm{d} x_{n-1}^2+\mathrm{d} y^2}{y^2}
\end{equation*}
and the distance formula 
\begin{equation*}
\cosh \eta(z',z)=1+\frac{||z'-z||^2}{2 y y'}
\end{equation*}
with $\eta:=\eta(O,z)$. The Laplace operator on $\mathbb{H}^n$ in cartesian coordinates is given by 
\begin{equation*}
\Delta_n=y^2 \left( \sum_{i=1}^{n-1} \frac{\partial^2}{\partial x_i^2}+\frac{\partial^2}{\partial y^2}\right)-(n-2)y \frac{\partial}{\partial y}
\end{equation*}
(for a proof see, for example, Chavel \cite{chavel} page 265). The Laplacian $\Delta_n$ in hyperbolic coordinates $(\eta, \bm{\alpha})=(\eta, \alpha_1, \dots, \alpha_{n-1})$, reads 
\begin{align} \label{radial}
\Delta_n=\frac{\partial^2}{\partial \eta^2}+\frac{n-1}{\tanh \eta} \frac{\partial}{\partial \eta} + \frac{1}{\sinh^2 \eta} \Delta_{S_{n-1}}
\end{align}
where $\Delta_{S_{n-1}}$ is the Laplace operator on the $(n-1)$-dimensional unit sphere (see, for example, Helgason \cite{helgason} page 158). 

In Lemma \ref{Acoth} we evaluate the hyperbolic Laplacian of the distance $\eta$ in $\mathbb{H}^n$. This result permits us, in Theorem \ref{Acoth2}, to determine the hyperbolic Laplacian of a smooth function $f(\eta)$. The statement of this result is given, for example, in Davies \cite{davies} page 117 without proof where $\cosh \rho$ must be replaced by $\coth \rho$.

\begin{lem} \label{Acoth}
For $z=(x,y)$ and $z'=(x',y')$ in $\mathbb{H}^n$ we have that the hyperbolic distance $\eta(z,z')$ is a solution of
\begin{equation*}
\Delta_n \eta(z,z')=\frac{n-1}{\tanh \eta(z,z')}.
\end{equation*}
\end{lem}
\Dim
Since $\coth( \operatorname{arcosh} (x))=\frac{x}{\sqrt{x^2-1}}$ and $\eta(z,z')= \operatorname{arcosh}\frac{||x-x'||^2+y^2+y'^2}{2 y y'},$
we have to prove that 
\begin{eqnarray*}
\Delta_n \eta(z,z')=(n-1) \frac{ ||x-x'||^2 +y^2+y'^2  }{\sqrt{ [\;||x-x'||^2+(y+y')^2\;]  \; [\;||x-x'||^2+(y-y')^2\;]  }}.
\end{eqnarray*}
In fact we have 
\begin{eqnarray} \label{derivataprima}
\frac{\partial}{\partial y} \eta(z,z')&=&-\frac{ ||x-x'||^2+y'^2-y^2   }{y\sqrt{ [\;||x-x'||^2+(y+y')^2\;]  \; [\;||x-x'||^2+(y-y')^2\;]  }},\\
\frac{\partial^2}{\partial y^2} \eta(z,z')&=&\frac{ ||x-x'||^2+y'^2-y^2   }{y^2\sqrt{ [\;||x-x'||^2+(y+y')^2\;]  \; [\;||x-x'||^2+(y-y')^2\;]  }}\nonumber\\
&&+\frac{ 4 ||x-x'||^2[||x-x'||^2+y'^2+y^2]   }{\sqrt{ [\;||x-x'||^2+(y+y')^2\;]^3\;   [\;||x-x'||^2+(y-y')^2\;]^3 }}. \nonumber
\end{eqnarray}
On the other side, for $i=1,\dots, n-1$, we have 
\begin{eqnarray} \label{derivataseconda}
\frac{\partial}{\partial x_i} \eta(z,z')&=&\frac{2(x_i-x'_i)}{\sqrt{ [\;||x-x'||^2+(y+y')^2\;] \;   [\;||x-x'||^2+(y-y')^2\;]  }},\\
\frac{\partial^2}{\partial x_i^2} \eta(z,z')&=&\frac{2}{\sqrt{ [\;||x-x'||^2+(y+y')^2\;]  \;  [\;||x-x'||^2+(y-y')^2\;]  }}\nonumber \\
&&-\frac{4(x_i-x'_i)^2[\;||x-x'||^2+y^2+y'^2\;]}{\sqrt{ [\;||x-x'||^2+(y+y')^2\;]^3 \;   [\;||x-x'||^2+(y-y')^2\;]^3  }},\nonumber \\
\sum_{i=1}^{n-1} \frac{\partial^2}{\partial x_i^2}\eta(z,z')&=&\frac{2(n-1)}{\sqrt{ [\;||x-x'||^2+(y+y')^2\;]   \; [\;||x-x'||^2+(y-y')^2\;]  }}\nonumber \\
&&-\frac{4||x-x'||^2[\;||x-x'||^2+y^2+y'^2\;]}{\sqrt{ [\;||x-x'||^2+(y+y')^2\;]^3  \;  [\;||x-x'||^2+(y-y')^2\;]^3  }}.\nonumber
\end{eqnarray}
So, finally, we obtain that 
\begin{eqnarray*}
\Delta_n \eta(z,z')&=&y^2 \left[ \sum_{i=1}^{n-1} \frac{\partial^2}{\partial x_i^2} \eta(z,z') +\frac{\partial^2}{\partial y^2} \eta(z,z')   \right]-(n-2) y  \frac{\partial}{\partial y}\eta(z,z')\\
&=&\frac{2(n-1)y^2}{\sqrt{ [\;||x-x'||^2+(y+y')^2\;] \;   [\;||x-x'||^2+(y-y')^2\;]  }}\\
&&+\frac{ ||x-x'||^2+y'^2-y^2   }{\sqrt{ [\;||x-x'||^2+(y+y')^2\;]  \; [\;||x-x'||^2+(y-y')^2\;]  }}\\
&&+\frac{ (n-2)[ \;||x-x'||^2+y'^2-y^2\;]   }{\sqrt{ [\;||x-x'||^2+(y+y')^2\;] \;  [\;||x-x'||^2+(y-y')^2\;]  }}\\
&=&(n-1) \frac{ ||x-x'||^2 +y'^2+y^2  }{\sqrt{ [\;||x-x'||^2+(y+y')^2\;] \;  [\;||x-x'||^2+(y-y')^2\;]  }}.
\end{eqnarray*}
\EndDim

\noindent In view of Lemma \ref{Acoth} the following theorem holds:  

\begin{teo} \label{Acoth2}
If $f$ is a smooth function on $\mathbb{R}$, it holds that  
\begin{equation*}
\Delta_n f(\eta)=f''(\eta)+\frac{n-1}{\tanh \eta} f'(\eta).
\end{equation*}
\end{teo}
\Dim
We have 
\begin{align*}
&\Delta_n f(\eta(z,z'))\\
&=y^2 \left[ \sum_{i=1}^{n-1} \frac{\partial^2}{\partial x_i^2} f(\eta(z,z')) +\frac{\partial^2}{\partial y^2} f (\eta(z,z'))   \right]-(n-2) y  \frac{\partial}{\partial y}f(\eta(z,z'))\\
&= y^2 \left[  \sum_{i=1}^{n-1}\left(\frac{\partial^2 f}{\partial \eta^2} \left(\frac{\partial \eta}{\partial x_i }\right)^2+\frac{\partial f}{\partial \eta} \frac{\partial^2 \eta}{\partial x_i^2}  \right)  + \frac{\partial^2 f}{\partial \eta^2} \left(\frac{\partial \eta}{\partial y }\right)^2+\frac{\partial f}{\partial \eta} \frac{\partial^2 \eta}{\partial y^2} \right] -(n-2)y \frac{\partial f}{\partial \eta} \frac{\partial \eta}{\partial y}\\
&=  \frac{\partial^2 f}{\partial \eta^2} y^2 \left[ \sum_{i=1}^{n-1}  \left(\frac{\partial \eta}{\partial x_i }\right)^2+ \left(\frac{\partial \eta}{\partial y }\right)^2 \right]+ \frac{\partial f}{\partial \eta} \left[ y^2 \left(   \sum_{i=1}^{n-1}   \frac{\partial^2 \eta}{\partial x_i^2 }+ \frac{\partial^2 \eta}{\partial y^2 } \right)-(n-2) y \frac{\partial \eta}{\partial y} \right]\\
&= \frac{\partial^2 f}{\partial \eta^2}\;y^2 \left[ \sum_{i=1}^{n-1}  \left(\frac{\partial \eta}{\partial x_i }\right)^2+ \left(\frac{\partial \eta}{\partial y }\right)^2   \right]+ \frac{\partial f}{\partial \eta} \Delta_n \eta\\
&= \frac{\partial^2 f}{\partial \eta^2} + \frac{\partial f}{\partial \eta} \Delta_n \eta,
\end{align*}
since, in view of formula (\ref{derivataprima}) and (\ref{derivataseconda}), it holds that 
\begin{eqnarray*}
\sum_{i=1}^{n-1}  \left(\frac{\partial \eta}{\partial x_i }\right)^2+ \left(\frac{\partial \eta}{\partial y }\right)^2&=&\frac{4 ||x-x'||^2 y^2+ [\;||x-x'||^2+y'^2-y^2\;]^2}{y^2[\;||x-x'||^2+(y+y')^2\;] \;  [\;||x-x'||^2+(y-y')^2\;]}=\frac{1}{y^2}.
\end{eqnarray*}
From Lemma \ref{Acoth} we obtain the final result 
$$\Delta_n f(\eta(z,z'))= \frac{\partial^2 f}{\partial \eta^2} + \frac{\partial f}{\partial \eta} \Delta_n \eta=\frac{\partial^2 f}{\partial \eta^2} + \frac{n-1}{\tanh \eta} \frac{\partial f}{\partial \eta}.
$$
\EndDim

We denote with $\{B_{\mathbb{H}^n}(t), t \ge 0\}$ the hyperbolic Brownian motion on $\mathbb{H}^n$ with starting point $z=(\eta,\bm{\alpha})\in \mathbb{H}^n$ where $\bm{\alpha}=(\alpha_1,\dots, \alpha_{n-1})\in [0,\pi]^{n-2} \times [0,2 \pi)$, and we assume that $z$ is inside the $n$-dimensional hyperbolic ball $U$ with hyperbolic radius $\bar{\eta}$ (see Figure \ref{sxx}). 

We are interested in obtaining the law of the position occupied by the hyperbolic Brownian motion on $\mathbb{H}^n$ when it hits the boundary $\partial U$ for the first time. 

Since the Laplace operator is invariant under rotations (see, for example, Helgason \cite{helgason} Proposition 2.4), without loss of generality we can assume that the starting point is $z=(\eta,\alpha_1,0, \dots,0)$ and the process hits the boundary of the ball $U$ at some point $\bar{z}=(\bar{\eta},\bar{\alpha}_1,0,\dots,0)$, where $\alpha_1-\bar{\alpha}_1$ is the angle between the vectors $z$ and $\bar{z}$. For a function on the $(n-1)$-dimensional unit sphere $S_{n-1}$ depending only on one angle ${\theta}$ we have 
\begin{align} \label{polar}
\Delta_{S_{n-1}}=\frac{1}{\sin^{n-2}\theta} \frac{\partial}{\partial \theta}\left(\sin^{n-2} \theta \frac{\partial}{\partial \theta} \right)=\frac{\partial^2}{\partial \theta^2}+\frac{n-2}{\tan \theta} \frac{\partial}{\partial \theta}.
\end{align}  
In view of (\ref{radial}) and (\ref{polar}) we have that the hitting distribution on $\partial U$ is obtained from the solution of the following Dirichlet problem 
\begin{equation} \label{2.34}
\begin{cases}
\left[ \frac{\partial^2}{\partial \eta^2}+\frac{n-1}{\tanh \eta} \frac{\partial}{\partial \eta}+\frac{1}{\sinh^2 \eta} \left( \frac{\partial^2}{\partial \alpha_1^2}+\frac{n-2}{\tan \alpha_1} \frac{\partial}{\partial \alpha_1} \right) \right] u(\eta,\alpha_1; \bar{\eta}, \bar{\alpha}_1)=0, &0<\eta < \bar{\eta}<\infty,\\
u(\bar{\eta},\alpha_1; \bar{\eta},\bar{\alpha}_1)=\delta(\alpha_1-\bar{\alpha}_1), &  \alpha_1-\bar{\alpha}_1 \in (0, \pi].
\end{cases}
\end{equation}

\begin{teo}
The solution to the Dirichlet problem (\ref{2.34}) is given by 
\begin{align} \label{petap}
&u(\eta, \alpha_1; \bar{\eta}, \bar{\alpha}_1) \nonumber \\
&\;=\frac{\Omega_{n-1}}{\Omega_n}   \sum_{k=0}^\infty \left(\frac{2k}{n-2}+1 \right)   \frac{  \tanh^k \frac{{\eta}} 2\; F\left(k,1-\frac n 2;k+\frac n 2;  \tanh^2 \frac{ {\eta}} 2 \right)}{ \tanh^k \frac{ \bar{\eta}} 2\; F\left(k,1-\frac n 2;k+\frac n 2;  \tanh^2 \frac{ \bar{\eta}} 2 \right)} C^{(\frac{n-2}{2})}_k(\cos( \alpha_1-\bar{\alpha}_1)) \sin^{n-2} (\alpha_1-\bar{\alpha}_1),
\end{align}
where $n>2$, $0<\eta < \bar{\eta}<\infty$ and $\alpha_1-\bar{\alpha}_1 \in (0, \pi]$.
\end{teo}
\Dim
As in Theorem \ref{theteo} our proof is based on the method of separation of variables. We assume that 
$$u(\eta,\alpha_1; \bar{\eta}, \bar{\alpha}_1)=\Theta(\alpha_1) E(\eta).$$
Since we have that 
\begin{align*}
 \Theta(\alpha_1) E''(\eta)+\Theta(\alpha_1) \frac{n-1}{\tanh \eta} E'(\eta)+\frac{E(\eta)}{\sinh^2 \eta} \left[ \Theta''(\alpha_1)+\frac{n-2}{\tan \alpha_1}\Theta'(\alpha_1)   \right]=0,
\end{align*}
there exists a constant $\mu^2$ such that 
\begin{equation} \label{firstsecond}
\begin{cases} 
\Theta''(\alpha_1)+(n-2) \cot \alpha_1 \; \Theta'(\alpha_1)+\mu^2 \Theta(\alpha_1)=0,  \\
\sinh^2 \eta \; E''(\eta) +(n-1) \cosh \eta\; \sinh \eta\;E'(\eta)-\mu^2 E(\eta)=0. 
\end{cases}
\end{equation}   
The first equation in (\ref{firstsecond}) can be reduced to the Gegenbauer equation. With the change of variable $\omega=\cos \alpha_1$ and for $\mu^2=k(k+n-2)$, we obtain
\begin{align} \label{pppu}
(1-\omega^2) G''(\omega)-(n-1) \omega G'(\omega)+k(k+n-2) G(\omega)=0.
\end{align} 
The Gegenbauer polynomials $C_k^{(\frac{n-2}{2})}(\omega)$ satisfy (\ref{pppu}) (see, for example, Polyanin and Zaitsev \cite{polyanin} S.2.11-4 or Helgason \cite{helgason} page 16) and this implies that 
\begin{align} \label{(3)}
\Theta(\alpha_1)=A C_k^{(\frac{n-2}{2})}(\cos \alpha_1).
\end{align}
We transform the second equation of (\ref{firstsecond}) 
\begin{align} \label{unouno}
\sinh^2 \eta E''(\eta) +(n-1) \cosh \eta \sinh \eta E'(\eta)-\mu^2E(\eta)=0
\end{align}
into a hypergeometric equation. The first step is based on the change of variable  $\zeta=\tanh \frac{\eta}{2}$. We have that 
\begin{align*}
\frac{\mathrm{d}}{\mathrm{d}\eta}&= \frac{1}{2 \cosh^2 \frac \eta 2} \frac{\mathrm{d}}{\mathrm{d}\zeta},\\
\frac{\mathrm{d}^2}{\mathrm{d}\eta^2}&= \frac{1}{4 \cosh^4 \frac \eta 2} \frac{\mathrm{d}^2}{\mathrm{d}\zeta^2}-\frac{\sinh\frac  \eta 2}{2 \cosh^3 \frac \eta 2} \frac{\mathrm{d}}{\mathrm{d}\zeta}.
\end{align*}
By taking into account that $\sinh \eta=2 \sinh \frac \eta 2 \cosh \frac \eta 2$ and that  $\cosh \eta=2 \cosh^2 \frac \eta 2-1$, equation (\ref{unouno}) becomes 
\begin{align*}
&\hspace{-1.5cm}4 \sinh^2 \frac{\eta}{2} \cosh^2 \frac \eta 2 \left[ \frac{1}{4 \cosh^4 \frac \eta 2}  \frac{\mathrm{d}^2}{\mathrm{d}\zeta^2}-\frac{\sinh \frac \eta 2}{2 \cosh^3 \frac \eta 2} \frac{\mathrm{d}}{\mathrm{d}\zeta}  \right] E(\zeta)\\
&+(n-1) 2 \sinh \frac \eta 2 \cosh \frac \eta 2 \left(2 \cosh^2 \frac \eta 2-1\right) \frac{1}{2 \cosh^2 \frac \eta 2} \frac{\mathrm{d}}{\mathrm{d}\zeta} E(\zeta)-\mu^2 E(\zeta)=0.
\end{align*}
And since
\begin{align*}
-\frac{2 \sinh^3 \frac  \eta 2}{\cosh \frac \eta 2}+(n-1) \frac{\sinh \frac \eta 2}{\cosh \frac \eta 2} \left( 2 \cosh^2 \frac \eta 2-1 \right)&=\tanh \frac \eta 2 \left[ -2 \sinh^2 \frac \eta 2+(n-1) \left( 2 \cosh^2 \frac \eta 2-1 \right)  \right]\\
&=\tanh \frac \eta 2 \left[  1+(n-2) \left( 2 \cosh^2 \frac \eta 2-1  \right) \right]\\
&=\tanh \frac \eta 2 \left[  1+(n-2) \frac{1+\tanh^2 \frac \eta 2}{1-\tanh^2 \frac \eta 2} \right],
\end{align*}
we can write equation (\ref{unouno}) as
\begin{align*} 
\tanh^2 \frac \eta 2 E''(\zeta)+\tanh \frac \eta 2 \left[ 1+(n-2) \frac{1+\tanh^2 \frac \eta 2}{1-\tanh^2 \frac \eta 2} \right] E'(\zeta)-\mu^2 E(\zeta)=0,
\end{align*}
that is
\begin{align} \label{unounouno}
\zeta^2 E''(\zeta)+\zeta \left[ 1+(n-2) \frac{1+\zeta^2}{1-\zeta^2} \right] E'(\zeta)-\mu^2 E(\zeta)=0.
\end{align}
We now assume that 
\begin{align*}
E(\zeta)=\zeta^k f(\zeta^2).
\end{align*}
Since 
\begin{align}
E'(\zeta)&=k \zeta^{k-1} f(\zeta^2)+2 \zeta^{k+1} f'(\zeta^2), \label{duo}\\
E''(\zeta)&=k(k-1) \zeta^{k-2} f(\zeta^2)+2 k \zeta^k f'(\zeta^2)+2(k+1)\zeta^k f'(\zeta^2)+4 \zeta^{k+2}f''(\zeta^2)\nonumber\\
&=k(k-1) \zeta^{k-2} f (\zeta^2)+2 (2 k+1) \zeta^k f'(\zeta^2)+4 \zeta^{k+2} f''(\zeta^2), \label{duoduo}
\end{align}
by replacing (\ref{duo}) and (\ref{duoduo}) into (\ref{unounouno}) (with $\mu^2=k(k+n-2)$) we have that 
\begin{align*}
\hspace{-0.5cm}&k(k-1) \zeta^k f(\zeta^2)+2(2k+1)\zeta^{k+2} f'(\zeta^2)+4 \zeta^{k+4} f''(\zeta^2)\\ 
&+\left[ 1+(n-2) \frac{1+\zeta^2}{1-\zeta^2} \right] \left[ k \zeta^k f(\zeta^2)+2 \zeta^{k+2} f'(\zeta^2) \right]-k(k+n-2) \zeta^k f(\zeta^2)=0
\end{align*}
and with obvious simplifications we have that 
\begin{align*}
&\hspace{-0.5cm}k(k-1)  f(\zeta^2)+2(2k+1)\zeta^{2} f'(\zeta^2)+4 \zeta^{4} f''(\zeta^2)\\ 
&+\left[ 1+(n-2) \frac{1+\zeta^2}{1-\zeta^2} \right] \left[ k f(\zeta^2)+2 \zeta^{2} f'(\zeta^2) \right]-k(k+n-2) f(\zeta^2)\\
&\hspace{-0.5cm}=4 \zeta^4 f''(\zeta^2)+2 \zeta^2\left[  2(k+1) +(n-2) \frac{1+\zeta^2}{1-\zeta^2}\right] f'(\zeta^2)+2 (n-2) k \frac{\zeta^2}{1-\zeta^2} f(\zeta^2)=0.
\end{align*}
After some additional manipulations we arrive at the following equation  
\begin{align} \label{2.444}
 \zeta^2 (1-\zeta^2) f''(\zeta^2)+\left[  k+\frac n 2 -  \left(k+2- \frac n 2\right) \zeta^2\right] f'(\zeta^2)+  k \left(\frac n 2 -1\right) f(\zeta^2)=0.
\end{align} 
Equation (\ref{2.444}) coincides with the hypergeometric equation $$t(1-t)f''(t)+[\gamma-(\alpha+\beta+1)t]f'(t)-\alpha \beta f(t)=0$$ for $t=\zeta^2$, $\alpha=k$, $\beta=1-\frac n 2$ and $\gamma=k+ \frac n 2$. In view of the position $E(\zeta)=\zeta^k f(\zeta^2)$ and $\zeta=\tanh \frac \eta 2$, we conclude that a solution to (\ref{unouno}) is given by 
\begin{align} \label{(4)}
E(\eta)= \tanh^k \frac \eta 2\; F\left(k,1-\frac n 2;k+\frac n 2;\tanh^2 \frac \eta 2  \right).
\end{align}
Equations (\ref{(3)}) and (\ref{(4)}) imply that 
\begin{align*}
u(\eta,\alpha_1; \bar{\eta}, \bar{\alpha}_1)&=\sum_{k=0}^\infty E_k(\eta) \Theta_k(\alpha_1) \\
&=\sum_{k=0}^\infty A_k  \tanh^k \frac \eta 2\; F\left(k,1-\frac n 2;k+\frac n 2; \tanh^2 \frac \eta 2 \right) C_k^{(\frac{n-2}{2})}(\cos \alpha_1).
\end{align*}
In order to determine the coefficients $A_k$ by applying the boundary conditions we have that  
\begin{align*}
u({\bar{ \eta}},\alpha_1; \bar{\eta}, \bar{\alpha}_1)=\delta(\alpha_1-\bar{\alpha}_1)
=\sum_{k=0}^\infty A_k   \tanh^k \frac{ \bar{\eta}} 2\;  F\left(k,1-\frac n 2;k+\frac n 2; \tanh^2 \frac{ \bar{\eta}} 2\right) C_k^{(\frac{n-2}{2})}(\cos \alpha_1).
\end{align*}
By multiplying both members by $C^{(\frac{n-2}{2})}_m(\cos \alpha_1)  \sin^{n-2}\alpha_1$  and then integrating we have that  
\begin{align*}
&\int_{0}^\pi \delta(\alpha_1-\bar{\alpha}_1) C^{(\frac{n-2}{2})}_m(\cos \alpha_1) \sin^{n-2}\alpha_1\; \mathrm{d}\alpha_1\\
&=\sum_{k=0}^\infty A_k  \tanh^k \frac{ \bar{\eta}} 2\; F\left(k,1-\frac n 2;k+\frac n 2; \tanh^2 \frac{ \bar{\eta}} 2 \right) \int_0^\pi C_k^{(\frac{n-2}{2})}(\cos \alpha_1) C^{(\frac{n-2}{2})}_m(\cos \alpha_1) \sin^{n-2}\alpha_1\; \mathrm{d}\alpha_1\\
&=A_m  \tanh^m \frac{ \bar{\eta}} 2\; F\left(m,1-\frac n 2;m+\frac n 2; \tanh^2 \frac{ \bar{\eta}} 2 \right) \frac{\pi 2^{3-n} \Gamma(m+n-2)}{m! \left( m+\frac{n-2}{2} \right) \Gamma\left(\frac{n-2}{2}  \right)^2  },
\end{align*} 
because the functions $C^{(n)}_k(x)$ form an orthogonal system on the interval $x \in (-1,1)$ (see Gradshteyn and Ryzhik \cite{G-R} formula 7.313). This implies that
\begin{align*}
A_m= \frac{C^{(\frac{n-2}{2})}_m(\cos \bar{\alpha}_1) \sin^{n-2} \bar{\alpha}_1}{    \tanh^m \frac{ \bar{\eta}} 2\;     F\left(m,1-\frac n 2;m+\frac n 2;  \tanh^2 \frac{ \bar{\eta}} 2\; \right)} \frac{m! \left( m+\frac{n-2}{2} \right) \Gamma\left(\frac{n-2}{2}  \right)^2  }{\pi 2^{3-n} \Gamma(m+n-2)}.
\end{align*}
We finally obtain 
\begin{align*}
&u(\eta,\alpha_1; \bar{\eta}, \bar{\alpha}_1)\\
&=\frac{\Gamma\left(\frac{n-2}{2}  \right)^2\sin^{n-2} \bar{\alpha}_1 }{2^{3-n} \pi}   \sum_{k=0}^\infty \frac{k! \left( k+\frac{n-2}{2} \right)   }{\Gamma(k+n-2)}  \frac{ \tanh^k \frac{ {\eta}} 2\; F\left(k,1-\frac n 2;k+\frac n 2;  \tanh^2 \frac{{\eta}} 2 \right)}{ \tanh^k \frac{ \bar{\eta}} 2\; F\left(k,1-\frac n 2;k+\frac n 2;  \tanh^2 \frac{ \bar{\eta}} 2 \right)} C^{(\frac{n-2}{2})}_k(\cos \bar{\alpha}_1) C_k^{(\frac{n-2}{2})}(\cos \alpha_1).
\end{align*}
By rotational invariance and since $C_k^{(\frac{n-2}{2})}(1)=\binom{n+k-3}{k}$, the last expression reduces to
\begin{align*}
&u(\eta, \alpha_1; \bar{\eta}, \bar{\alpha}_1)\\
&\;=\frac{\Gamma\left(\frac{n-2}{2}  \right)^2\sin^{n-2} (\alpha_1-\bar{\alpha}_1) }{2^{3-n} (n-3)! \pi}   \sum_{k=0}^\infty \left(k+\frac{n-2}{2} \right)  \frac{  \tanh^k \frac{ {\eta}} 2\; F\left(k,1-\frac n 2;k+\frac n 2;  \tanh^2 \frac{ {\eta}} 2 \right)}{ \tanh^k \frac{ \bar{\eta}} 2\; F\left(k,1-\frac n 2;k+\frac n 2;  \tanh^2 \frac{ \bar{\eta}} 2 \right)} C^{(\frac{n-2}{2})}_k(\cos( \alpha_1-\bar{\alpha}_1)) .
\end{align*}
We arrive at formula (\ref{petap}) by observing that 
\begin{align*}
\frac{\Gamma\left(\frac{n-2}{2}  \right)^2 }{2^{3-n} (n-3)! \pi}=\frac{2}{n-2} \frac{\Omega_{n-1}}{\Omega_n}
\end{align*}
where $\Omega_n=\frac{2 \pi^{n/2}}{\Gamma(n/2)}$ is the surface area of the $n$-dimensional Euclidean unit sphere. In fact, since $2^{2-n} \sqrt \pi \Gamma(n-1)=\Gamma(\frac n 2) \Gamma(\frac n 2-\frac 1 2)$, we have 
\begin{align*}
\frac{\Gamma\left(\frac{n-2}{2}  \right)^2 }{2^{3-n} (n-3)! \pi}&=\frac{1}{2^{3-n} (n-3)! \pi} \frac{\left( \frac n 2-1 \right)^2}{\left( \frac n 2-1 \right)^2} \Gamma\left(\frac{n-2}{2}  \right)^2=\frac{1}{2^{3-n} (n-3)! \pi} \frac{\Gamma\left(\frac{n}{2}  \right)^2}{\left( \frac n 2-1 \right)^2}  \\
&=  \frac{2}{n-2} \frac{\Gamma\left( \frac{n}{2} \right)}{\sqrt \pi} \frac{\Gamma\left( \frac{n}{2} \right)}{2^{2-n} \sqrt \pi \Gamma(n-1)}=  \frac{2}{n-2} \frac{\Gamma\left( \frac{n}{2} \right)}{\sqrt \pi} \frac{1}{\Gamma\left(\frac n 2 -\frac 1 2 \right)}\\
&= \frac{2}{n-2} \frac{\Omega_{n-1}}{\Omega_n}
\end{align*}
and this concludes the proof of the theorem. 
\EndDim

\begin{rem} \label{theremark}
We note that for small values of $\eta$ and $\bar{\eta}$ we obtain the Euclidean Poisson kernel. In fact, since  $\tanh \frac \eta 2 \sim \frac \eta 2$, $C_1^{(n)}(t)=2 n t$, $C_0^{(n)}(t)=1$ and $kC_k^{(n)}(t)=2 n [tC_{k-1}^{(n+1)}(t)-C_{k-2}^{(n+1)}(t)]$, we have that 
\begin{align*}
&\hspace{-0.25cm}\sum_{k=0}^\infty \left(\frac{2k}{n-2}+1 \right)   \frac{  \tanh^k \frac{{\eta}} 2\; F\left(k,1-\frac n 2;k+\frac n 2;  \tanh^2 \frac{ {\eta}} 2 \right)}{ \tanh^k \frac{ \bar{\eta}} 2\; F\left(k,1-\frac n 2;k+\frac n 2;  \tanh^2 \frac{ \bar{\eta}} 2 \right)} C^{(\frac{n-2}{2})}_k(\cos( \alpha_1-\bar{\alpha}_1))\\
&  \sim   \sum_{k=0}^\infty \left(\frac{2k}{n-2}+1 \right) \left(\frac{\eta}{\bar{\eta}}\right)^k  C^{(\frac{n-2}{2})}_k(\cos (\alpha_1-\bar{\alpha}_1)) \\
&= \frac{2}{n-2}  \left[  \sum_{k=2}^\infty k \left(\frac{\eta}{\bar{\eta}}\right)^k  C^{(\frac{n-2}{2})}_k(\cos (\alpha_1-\bar{\alpha}_1)) + (n-2) \frac{\eta}{\bar{\eta}}  \cos (\alpha_1-\bar{\alpha}_1)  \right.\\
&\;\;\;\;+\left.\frac{n-2}{2} \sum_{k=0}^\infty  \left(\frac{\eta}{\bar{\eta}}\right)^k C^{(\frac{n-2}{2})}_k(\cos (\alpha_1-\bar{\alpha}_1))\right]\\
&=2 \left[ \cos (\alpha_1-\bar{\alpha}_1) \sum_{k=2}^\infty  \left(\frac{\eta}{\bar{\eta}}\right)^k  C^{(\frac{n}{2})}_{k-1}(\cos (\alpha_1-\bar{\alpha}_1))    \right.\\
&\;\;\;\;\left.-\sum_{k=2}^\infty   \left(\frac{\eta}{\bar{\eta}}\right)^k C^{(\frac{n}{2})}_{k-2}(\cos (\alpha_1-\bar{\alpha}_1))  +\frac{\eta}{\bar{\eta}} \cos (\alpha_1-\bar{\alpha}_1)+  \frac{1}{2} \sum_{k=0}^\infty  \left(\frac{\eta}{\bar{\eta}}\right)^k C^{(\frac{n-2}{2})}_k(\cos (\alpha_1-\bar{\alpha}_1)) \right]\\
&=2  \left[ \frac{\eta}{\bar{\eta}}  \cos(\alpha_1-\bar{\alpha}_1) \sum_{k=0}^\infty  \left(\frac{\eta}{\bar{\eta}}\right)^k  C^{(\frac{n}{2})}_{k}(\cos (\alpha_1-\bar{\alpha}_1))    \right.\\
&\;\;\;\;\left.-  \left(\frac{\eta}{\bar{\eta}}\right)^2 \sum_{k=0}^\infty   \left(\frac{\eta}{\bar{\eta}}\right)^k C^{(\frac{n}{2})}_{k}(\cos (\alpha_1-\bar{\alpha}_1)) +  \frac{1}{2} \sum_{k=0}^\infty  \left(\frac{\eta}{\bar{\eta}}\right)^k C^{(\frac{n-2}{2})}_k(\cos (\alpha_1-\bar{\alpha}_1)) \right]\\
&=2 \left[   \left( 1- \frac{2 \eta}{\bar{\eta}} \cos(\alpha_1-\bar{\alpha}_1)+   \frac{\eta^2}{\bar{\eta}^2} \right)^{-\frac n 2}  \left( \frac{\eta}{\bar{\eta}}  \cos(\alpha_1-\bar{\alpha}_1) -  \frac{\eta^2}{\bar{\eta}^2} \right) \right.\left.+  \frac{1}{2} \left( 1- \frac{2 \eta}{\bar{\eta}} \cos(\alpha_1-\bar{\alpha}_1)+   \frac{\eta^2}{\bar{\eta}^2} \right)^{-\frac{n-2} 2} \right]\\
&=2  \left( 1- \frac{2 \eta}{\bar{\eta}} \cos(\alpha_1-\bar{\alpha}_1)+   \frac{\eta^2}{\bar{\eta}^2} \right)^{-\frac n 2}    \left[  \frac{\eta}{\bar{\eta}}  \cos(\alpha_1-\bar{\alpha}_1) -  \frac{\eta^2}{\bar{\eta}^2}  +  \frac{1}{2} \left( 1- \frac{2 \eta}{\bar{\eta}} \cos(\alpha_1-\bar{\alpha}_1)+   \frac{\eta^2}{\bar{\eta}^2} \right) \right]\\
&=   \left( 1- \frac{2 \eta}{\bar{\eta}} \cos(\alpha_1-\bar{\alpha}_1)+   \frac{\eta^2}{\bar{\eta}^2} \right)^{-\frac n 2}    \left( 1-  \frac{\eta^2}{\bar{\eta}^2}  \right)\\
&=   \frac{ 1-  \frac{\eta^2}{\bar{\eta}^2} }{ \left( 1- \frac{2 \eta}{\bar{\eta}} \cos(\alpha_1-\bar{\alpha}_1)+   \frac{\eta^2}{\bar{\eta}^2} \right)^{\frac n 2} }.
\end{align*}
\end{rem}

\begin{rem}
The kernel (\ref{petap}) represents the marginal, with respect to $\bar{\alpha}_2, \dots,\bar{\alpha}_{n-1}$, of the distribution of the position occupied by the hyperbolic Brownian motion $\{B_{\mathbb{H}^n}(t), t \ge 0\}$ starting from $z=(\eta,\bm{\alpha}) \in \mathbb{H}^n$ when it hits for the first time the boundary $\partial U$ of the $n$-dimensional hyperbolic hypersphere of radius $\bar{\eta}$.  For $z=(\eta,{\bf 0})$, such distribution is given by 
\begin{align} \label{peta2}
\mathbb{P}_z\{B_{\mathbb{H}^n}(T_{\bar{\eta}}) \in \mathrm{d}\bar{\bm{\alpha}} \}=\sum_{k=0}^\infty \left(\frac{2k}{n-2}+1 \right)   \frac{\tanh^k \frac \eta 2 \; F\left(k,1-\frac n 2;k+\frac n 2; \tanh^2 \frac \eta 2 \right)}{\tanh^k \frac{\bar{\eta}} 2 \; F\left(k,1-\frac n 2;k+\frac n 2;\tanh^2 \frac{\bar{\eta}}2 \right)} C^{(\frac{n-2}{2})}_k(\cos \bar{\alpha}_1) f(\bar{\bm{\alpha}})  \mathrm{d}\bar{\bm{\alpha}} ,
\end{align} 
where $n>2$, $\eta< \bar{\eta}$ , $\bar{\alpha}_1 \in [0,\pi)$ is the angle between $z$ and $\bar{z}$, and
$$f(\bar{\bm{\alpha}})=\frac{1}{\Omega_n} \sin^{n-2} \bar{\alpha}_1 \sin^{n-3} \bar{\alpha}_2\dots \sin \bar{\alpha}_{n-2}$$
is the uniform density on $S_{n-1}$.
\end{rem}

\begin{rem}
We observe that (\ref{peta2}) is a proper probability law. In fact:
\begin{itemize}
\item  The non negativity is due to the non negativity of solutions of Dirichlet problems with non-negative boundary conditions.
\item It integrates to one, in fact
\end{itemize}
\begin{align*}
&\int_0^\pi   \dots \int_0^\pi  \int_0^{2 \pi}  \mathbb{P}_z\{B_{\mathbb{H}^n}(T_{\bar{\eta}}) \in \mathrm{d}\bar{\bm{\alpha}} \} \\
&=\frac{\Omega_{n-1}}{\Omega_n} \sum_{k=0}^\infty \left(\frac{2k}{n-2}+1 \right)  \frac{\tanh^k \frac  \eta 2 \; F\left(k,1-\frac n 2;k+\frac n 2; \tanh^2 \frac \eta 2 \right)}{\tanh^k\frac{\bar{\eta}}2\; F\left(k,1-\frac n 2;k+\frac n 2;\tanh^2 \frac{\bar{\eta}}2 \right)} \int_0^\pi  \; C^{(\frac{n-2}{2})}_k(\cos \bar{\alpha}_1)\; \sin^{n-2} \bar{\alpha}_1 \; \mathrm{d}\bar{\alpha}_1 \\
&=\frac{\Omega_{n-1}}{\Omega_n}   \frac{ F\left(0,1-\frac n 2;\frac n 2; \tanh^2 \frac \eta 2 \right)}{ F\left(0,1-\frac n 2;\frac n 2; \tanh^2 \frac{\bar{\eta}}2 \right)}  \int_0^\pi  \; C^{(\frac{n-2}{2})}_0(\cos \bar{\alpha}_1)\; \sin^{n-2} \bar{\alpha}_1 \; \mathrm{d}\bar{\alpha}_1 \\
&=\frac{\Omega_{n-1}}{\Omega_n}    \int_0^\pi \sin^{n-2} \bar{\alpha}_1  \; \mathrm{d}\bar{\alpha}_1 =1,
\end{align*}
since, if $k>0$, we have $$\int_0^\pi \; C^{(n)}_k(\cos \theta) \; \sin^{n-2} \theta \; \mathrm{d}\theta=0$$  (see Gradshteyn and Ryzhik \cite{G-R} formula 7.311.1) and $F(0,\beta;\gamma;z)=1$, $C_0^{n}(x)=1$,
$$\int_0^\pi \sin^{n-2} \theta \;  \mathrm{d}\theta=B\left(\frac 1 2, \frac{n-1}{2}\right).$$
\end{rem}

\begin{rem}
We also note that: 
\begin{itemize}
\item For $\eta \to 0$ (i.e. when the starting point is the center of the hyperbolic hypersphere) formula (\ref{peta2}) becomes the uniform distribution on $S_{n-1}$ as expected.
\item For $\bar{\eta} \to \infty$, since $\tanh \frac{\bar{\eta}}{2} \to 1$ and $F(\alpha,\beta;\gamma;1)=\frac{\Gamma(\gamma) \Gamma(\gamma-\alpha-\beta)}{\Gamma(\gamma-\alpha) \Gamma(\gamma-\beta)}$ if $\gamma>\alpha+\beta$ (see Gradshteyn and Ryzhik \cite{G-R} formula 9.122.1), we have that 
\end{itemize}
\begin{align*} 
&\lim_{\bar{\eta} \to \infty}\mathbb{P}_z\{B_{\mathbb{H}^n}(T_{\bar{\eta}}) \in \mathrm{d}\bar{\bm{\alpha}} \}\\
&\;\;=\sum_{k=0}^\infty \left(\frac{2k}{n-2}+1 \right) \frac{\Gamma(k+n-1)}{\Gamma(k+\frac n 2)} \; \tanh^k \frac \eta 2 \;\;  F\left(k,1-\frac n 2;k+\frac n 2; \tanh^2 \frac \eta 2 \right) C^{(\frac{n-2}{2})}_k(\cos \bar{\alpha}_1) f(\bar{\bm{\alpha}})  \mathrm{d}\bar{\bm{\alpha}}.
\end{align*} 
\end{rem}
\begin{rem}
 Byczkowski et al. in \cite{byczkowski} provide an integral formula for the hyperbolic Poisson kernel of  the half-space $H_a=\{(x,y)\in \mathbb{H}^n: y>a\}$ for $n>2$, $a>0$, and show that for $a \to 0$, it converges to the Cauchy-type distribution
 $$\frac{\Gamma(n-1)}{\pi^{\frac{n-1}{2}}\Gamma(\frac{n-1}{2})} \left( \frac{y}{y^2+|x|^2} \right)^{n-1}.$$
\end{rem}
 
\section{Exit probabilities from a hyperbolic annulus in $\mathbb{H}^n$}

\subsection{Two dimensional case}

Suppose the hyperbolic Brownian motion $\{B_{\mathbb{H}^2}(t), t\ge 0\}$ starts at  $z=(\eta,\alpha) \in \mathbb{H}^2$ inside the hyperbolic annulus $A$ with radii $0<\eta_1<\eta_2<\infty$
$$A=\{(\eta,\alpha): \eta_1<\eta<\eta_2\}$$
(see Figure \ref{autobus}). We define the hitting times
$$T_{\eta_i}=\inf \{t>0: \eta(O,B_{\mathbb{H}^2}(t))=\eta_i\},\hspace{1cm}  i=1,2,$$
and $T=T_{\eta_1} \wedge T_{\eta_2}$. In the next theorem we evaluate the exit probabilities
$\mathbb{P}_z\{T_{\eta_1}<T_{\eta_2}\}$. Since these are given in terms of harmonic functions on the annulus $A$, they are closely related to the Dirichlet problem. 

\begin{teo} \label{unabellalabel}
Let $\{B_{\mathbb{H}^2}(t):t\ge0\}$ be a hyperbolic Brownian motion starting at $z=(\eta,\alpha) \in A$. The following result holds true
\begin{equation} \label{uku}
\mathbb{P}_z\{T_{\eta_1}<T_{\eta_2}\}=\frac{ \log   \tanh \frac{\eta_2}{2} - \log   \tanh \frac{\eta}{2} }{ \log   \tanh \frac{\eta_2}{2} - \log   \tanh \frac{\eta_1}{2} }, \hspace{1cm} \eta_1<\eta<\eta_2.
\end{equation}
\end{teo}
\Dim 
Since the probability in (\ref{uku}) is spherically symmetric we are lead to study the solution $v: (\eta_1,\eta_2) \to \mathbb{R}$ to the Laplace equation involving only the radial part:
$$\left[ \frac{\partial^2}{\partial \eta^2}+ \frac{1}{\tanh \eta} \frac{\partial}{\partial \eta} \right] v(\eta)=0$$ 
subjected to the boundary conditions $v(\eta_1)=1$ and $v(\eta_2)=0$.
With the change of variable $w=\cosh \eta$ we immediately get 
\begin{equation} \label{xwxw}
(1-w^2) K''(w)-2 w K'(w)=0,
\end{equation}
whose general solution is 
\begin{equation*} 
K(w)=C_1+C_2 \log \left| \frac{w-1}{w+1}  \right|
\end{equation*}
(see, for example, Polyanin and Zaitsev  \cite{polyanin}, Section 2.1.2 Formula 233 for $a=1$, $b=-1$, $\lambda=0$ and $\mu=0$). It follows that 
\begin{equation} \label{verme}
v(\eta)=C_1+C_2 \log \left( \frac{\cosh \eta-1}{\cosh \eta+1}  \right)=C_1+C_2 \log \tanh \frac{\eta} 2.
\end{equation}
By imposing the boundary conditions we get 
\begin{eqnarray*}
\mathbb{P}_z\{T_{\eta_1}<T_{\eta_2}\}=\frac{v(\eta_2)-v(\eta)}{v(\eta_2)-v(\eta_1)}=\frac{ \log  \tanh \frac{\eta_2} 2- \log \tanh \frac{\eta} 2}{ \log  \tanh \frac{\eta_2} 2- \log  \tanh \frac{\eta_1} 2}.
\end{eqnarray*}
\EndDim

Starting from (\ref{uku}) and letting $\eta_2$ go to infinity we have that Theorem \ref{unabellalabel} leads to the following corollary.

\begin{cor}
For any $z=(\eta,\alpha)$ outside the hyperbolic disc of radius $\eta_1$ and center $O$, we have 
\begin{equation} \label{gkg}
 \mathbb{P}_z\{T_{\eta_1}<\infty\}=\frac{ \log \left( \frac{\cosh \eta-1}{\cosh \eta+1}  \right)}{ \log \left( \frac{\cosh \eta_1-1}{\cosh \eta_1+1}  \right)}=\frac{\log \tanh \frac{\eta}{2}}{\log \tanh \frac{\eta_1}{2}},
  \hspace{1cm} \eta_1<\eta
. \end{equation}
\end{cor} 

It is possible to show with simple computations that the functions in (\ref{uku}) and (\ref{gkg}) are genuine probabilities since they vary in $(0,1)$.

 \subsection{Multidimensional case}

It is possible to generalize the exit probabilities from a hyperbolic annulus to the case of the $n$-th dimensional hyperbolic Brownian motion.

In order to evaluate the exit probabilities from the hyperbolic annulus $A$ in $\mathbb{H}^n$, with hyperbolic radii $\eta_1$ and $\eta_2$ with $\eta_1<\eta_2$, we are interested in obtaining a solution $v_n:(\eta_1,\eta_2) \to \mathbb{R}$ to the radial part of the hyperbolic Laplace equation in $\mathbb{H}^n$. We have proved in Lemma \ref{Acoth} and Theorem \ref{Acoth2} that it is equivalent to solve 
\begin{equation} \label{raddial} \left[ \frac{\mathrm{d}^2}{\mathrm{d}\eta^2}+ \frac{n-1}{\tanh \eta}  \frac{\mathrm{d}}{\mathrm{d}\eta} \right]v_n(\eta)=0.\end{equation}
In what follows we will assume that 
\begin{equation*}
c(n,0)=1, \hspace{1cm} c(n,k)= \frac{(n-3)(n-5)\cdots(n-2 k-1)}{(n-2)(n-4) \cdots (n-2k-2)}, \hspace{0.7cm}  k=1, \dots \frac{n-3}{2}.
\end{equation*}

\begin{teo}
For a hyperbolic Brownian motion $\{B_{\mathbb{H}^n}(t): t \ge 0\}$ started at $z=(\eta, \alpha) \in A$, we have that \\

\noindent For $n=3,5,7, \dots$ 
\begin{equation} \label{copti}
\mathbb{P}_z\{T_{\eta_1}<T_{\eta_2}\}=\frac{\sum_{k=0}^{\frac{n-3}{2} } (-1)^{k-1}  c(n,k)   \left[ \frac{\cosh \eta_2}{\sinh^{n-2k-2}\eta_2} - \frac{\cosh \eta}{\sinh^{n-2k-2}\eta} \right]}{\sum_{k=0}^{\frac{n-3}{2} } (-1)^{k-1} c(n,k)  \left[ \frac{\cosh \eta_2}{\sinh^{n-2k-2}\eta_2} - \frac{\cosh \eta_1}{\sinh^{n-2k-2}\eta_1} \right]}. 
\end{equation}

\noindent For $n=4,6,8, \dots$ 
\begin{equation*}
\mathbb{P}_z\{T_{\eta_1}<T_{\eta_2}\}=\frac{\sum_{k=0}^{\frac{n-4}{2} } (-1)^{k-1} c(n,k)   \left[ \frac{\cosh \eta_2}{\sinh^{n-2k-2}\eta_2} - \frac{\cosh \eta}{\sinh^{n-2k-2}\eta} \right] +(-1)^{\frac{n-2}{2}} \frac{(n-3)!!}{(n-2)!!} \log \frac{\tanh \frac{\eta_2}{2}}{\tanh \frac{\eta}{2}}}{\sum_{k=0}^{\frac{n-4}{2} } (-1)^{k-1} c(n,k)   \left[ \frac{\cosh \eta_2}{\sinh^{n-2k-2}\eta_2} - \frac{\cosh \eta_1}{\sinh^{n-2k-2}\eta_1} \right]+(-1)^{\frac{n-2}{2}} \frac{(n-3)!!}{(n-2)!!} \log \frac{\tanh \frac{\eta_2}{2}}{\tanh \frac{\eta_1}{2}}}.
\end{equation*}
\end{teo}
\Dim
The general solution to equation (\ref{raddial})
is given by 
\begin{equation*}
v_n(\eta)=C_1+C_2 \int \frac{1}{\sinh^{n-1}\eta} \mathrm{d}\eta.
\end{equation*}
For $n=2m+1$, $m=1,2,\dots$ we have  
\begin{align} \label{ri1}
v_n(\eta)&=C_1+C_2 \int \frac{1}{\sinh^{2 m} \eta} \; \mathrm{d}\eta \nonumber  \\
&=C_1+ C_2 \frac{\cosh \eta}{2m-1} \left[ -\frac{1}{\sinh^{2m-1}\eta}+\sum_{k=1}^{m-1} (-1)^{k-1} \frac{2^k(m-1)(m-2)\cdots(m- k)}{(2m-3)(2m-5) \cdots (2m-2k-1)}\frac{1}{\sinh^{2m-2k-1}\eta}  \right] \nonumber\\
&=C_1+C_2 \sum_{k=0}^{\frac{n-3}{2}} (-1)^{k-1} C(n,k) \frac{\cosh \eta}{\sinh^{n-2k-2}\eta }
\end{align}
(see Gradshteyn and Ryzhik \cite{G-R} formula 2.416.2).\\

\noindent For $n=2m+2$, $m=1,2,\dots$ we have
\begin{align} \label{ri2}
v_n(\eta)&=C_1+C_2 \int \frac{1}{\sinh^{2 m+1} \eta} \; \mathrm{d}\eta \nonumber \\
&=C_1+ C_2\frac{\cosh \eta}{2m} \left[ -\frac{1}{\sinh^{2m}\eta}+\sum_{k=1}^{m-1} (-1)^{k-1} \frac{(2m-1)(2m-3)\cdots(2m-2 k+1)}{2^k(m-1)(m-2) \cdots (m-k)}\frac{1}{\sinh^{2m-2k}\eta} \right] \nonumber \\
 & \hspace{0.4cm} +C_2  (-1)^m \frac{(2m-1)!!}{(2m)!!} \log \tanh \frac \eta 2 \nonumber \\
 &= C_1+ C_2\left[ \sum_{k=0}^{\frac{n-4}{2}} (-1)^{k-1} C(n,k) \frac{\cosh \eta}{\sinh^{n-2k-2}\eta}  +  (-1)^{\frac{n-2}{2}} \frac{(n-3)!!}{(n-2)!!} \log \tanh \frac \eta 2\right]
\end{align}
(see Gradshteyn and Ryzhik \cite{G-R} formula  2.416.3). With computations analogous to those performed in the two dimensional case, we obtain the statement. 
\EndDim

\noindent From this it follows immediately that:  

\begin{cor}
For $z=(\eta,\alpha)$ outside the hyperbolic  ball in $\mathbb{H}^n$ with radius $\eta_1$ and center in $O$, we have that \\

\noindent For $n=3,5,7, \dots$ 
\begin{equation*}
\mathbb{P}_z\{T_{\eta_1}<\infty\}=\frac{\sum_{k=0}^{\frac{n-5}{2} } (-1)^{k} c(n,k) \frac{\cosh \eta}{\sinh^{n-2k-2}\eta} +(-1)^{\frac{n-5}{2}}   \frac{(n-3)!!}{(n-4)!!} \left[ 1-\frac{\cosh \eta}{\sinh \eta} \right] }{\sum_{k=0}^{\frac{n-5}{2} } (-1)^{k} c(n,k) \frac{\cosh \eta_1}{\sinh^{n-2k-2}\eta_1} + (-1)^{\frac{n-5}{2}}   \frac{(n-3)!!}{(n-4)!!} \left[ 1-\frac{\cosh \eta_1}{\sinh \eta_1} \right]}, \hspace{1cm}\eta_1<\eta.
\end{equation*}

\noindent For $n=4,6,8, \dots$ 
\begin{equation*}
\mathbb{P}_z\{T_{\eta_1}<\infty\}=\frac{\sum_{k=0}^{\frac{n-4}{2} } (-1)^{k} c(n,k) \frac{\cosh \eta}{\sinh^{n-2k-2}\eta} +(-1)^{\frac n 2}   \frac{(n-3)!!}{(n-2)!!}  \log \tanh \frac {\eta}{2} }{\sum_{k=0}^{\frac{n-4}{2} } (-1)^{k} c(n,k)\frac{\cosh \eta_1}{\sinh^{n-2k-2}\eta_1} + (-1)^{\frac n 2} \frac{(n-3)!!}{(n-2)!!}  \log \tanh \frac {\eta_1}{2} }, \hspace{1cm}\eta_1<\eta.
\end{equation*}
\end{cor}

\begin{rem}
For the space $\mathbb{H}^3$ formula (\ref{copti}) takes the simple form
$$\mathbb{P}_z\{T_{\eta_1}<T_{\eta_2}\}=\frac{\coth \eta_2-\coth \eta}{\coth \eta_2-\coth \eta_1},\hspace{1cm}\eta_1<\eta<\eta_2,$$
and for $\eta_2 \to \infty$ yields 
$$\mathbb{P}_z\{T_{\eta_1}<\infty\}=\frac{1-\coth \eta}{1-\coth \eta_1}<1.$$
This shows that there is a positive probability that the hyperbolic Brownian motion never hits the ball of radius $\eta_1$.
\end{rem}

\begin{rem}
We note that for small values of $\eta$ we have $\frac{\cosh \eta}{\sinh^p \eta}\sim \frac{1}{\eta^p}$ and $\log \tanh \frac \eta 2\sim \log \eta $. From (\ref{verme}), (\ref{ri1}) and (\ref{ri2}) it follows that 
\begin{equation*}
v_n(\eta) \sim 
\begin{cases}
C_1+C_2 \log \eta, &\mathrm{if\;} n=2,\\
C_1+C_2 \eta^{2-n}, &\mathrm{if\;} n=3,4,5\dots
\end{cases}
\end{equation*}  
This means that, for sufficiently small domains, we obtain  the exit probabilities of Euclidean Brownian motion from an annulus:   
\begin{align} \label{sodac}
\mathbb{P}_z\{T_{\eta_1}<T_{\eta_2}\} \sim 
\begin{cases}
\vspace{0.2cm}
\frac{\log \eta_2 -\log \eta}{\log \eta_2-\log \eta_1}, &\mathrm{if\;} n=2,\\
\frac{ \eta_2^{2-n}-\eta^{2-n}}{\eta_2^{2-n}-\eta_1^{2-n}}, &\mathrm{if\;} n=3,4,5\dots
\end{cases}
\end{align}  
\end{rem}

\begin{rem}
It is important to note that for a planar hyperbolic Bownian motion the probability that the process goes to infinity before hitting the hyperbolic circle of radius $\eta_1$ is strictly less then one
$$ \mathbb{P}_z\{T_{\eta_1}<\infty\}=\frac{\log \tanh \frac{\eta}{2}}{\log \tanh \frac{\eta_1}{2}}<1$$
while it is well known, see (\ref{sodac}), that for a planar Euclidean Brownian motion it holds that 
$$\mathbb{P}_z\{T_{\eta_1}<\infty\}=1.$$ Hyperbolic Brownian motion is, in fact, transient for every dimension $n \ge 2$ as stated in Grigor'yan \cite{grigoryan} Proposition 3.2.   
\end{rem}

\section{Hitting distribution on a hyperbolic circle in $\mathbb{D}^2$}

It is possible to obtain analogous results by considering a different model of the hyperbolic plane. In particular in this section we consider the Poincar\'e disc model $\mathbb{D}^2$ instead of the half-plane model $\mathbb{H}^2$. The half-plane $\mathbb{H}^2$ can be mapped onto the disc $\mathbb{D}^2=\{(r,\theta): r \in [0,1), \theta \in (-\pi,\pi]\}$ by means of the conformal mapping $f: \mathbb{H}^2 \to \mathbb{D}^2$ such that
\begin{align}\label{xxx}f(z)=\frac{iz+1}{z+i}.\end{align}
The $x$-axis of $\mathbb{H}^2$ is mapped onto $\partial \mathbb{D}^2$ while the origin $O=(0,1)$ of $\mathbb{H}^2$ is mapped into the origin $O=(0,0)$ of $\mathbb{D}^2$. An arbitrary point $z=(x,y) \in \mathbb{H}^2$ is mapped into a point $Q=(r,\theta) \in \mathbb{D}^2$ such that 
\begin{equation}  \label{huhu}
\hspace{1cm}
\begin{cases}
x=\frac{2 r \cos \theta}{1+r^2-2 r \sin \theta},   \\
y=\frac{1-r^2}{1+r^2-2 r \sin \theta}
\end{cases}
\end{equation}
(for details see Lao and Orsingher \cite{lao}). In view of (\ref{hjhj}) and (\ref{huhu}) we have
\begin{equation*}
\frac{x}{y}=\sinh \eta \cos \alpha=\frac{2 r \cos \theta}{1-r^2}.
\end{equation*}
Since we have that
\begin{eqnarray*}
\cos \alpha=\cos \theta \hspace{1cm} \mathrm{and} \hspace{1cm} \sinh \eta=\frac{2 r}{1-r^2},
\end{eqnarray*}
for $\theta, \alpha \in (-\pi, \pi]$, we easily arrive at
\begin{eqnarray} \label{rr}
\begin{cases}
r=\frac{\cosh \eta-1}{\sinh \eta}=\sqrt{\frac{\cosh \eta-1}{\cosh \eta+1}}=\tanh \frac \eta 2,\\
\theta=\alpha.
\end{cases}
\end{eqnarray}
The hyperbolic metric and the distance formula in $\mathbb{D}^2$ become 
$$\mathrm{d}s^2=\frac{4}{(1-r^2)^2} \mathrm{d}r^2, \hspace{1cm}\mathrm{d}(O,Q)=\log \frac{1+r}{1-r}.$$
By means of (\ref{huhu}) the hyperbolic Laplacian in (\ref{1.5}) is converted into 
$$(1-r^2)^2 \left[ \frac{1}{r} \frac{\partial}{\partial r} \left( r \frac{\partial}{\partial r} \right)+\frac{1}{r^2} \frac{\partial^2}{\partial \theta^2} \right],$$
and the Dirichlet  problem for the hyperbolic disc $U=\{(r,\theta): r<\bar{r}\}$ in $\mathbb{D}^2$ reads 
\begin{equation} \label{gfgf}
\begin{cases}
(1-r^2)^2 \left[ \frac{1}{r} \frac{\partial}{\partial r} \left( r \frac{\partial}{\partial r} \right)+\frac{1}{r^2} \frac{\partial^2}{\partial \theta^2} \right] u(r, \theta;\bar{r}, \bar{\theta})=0, &0<r < \bar{r}<1, \\
u( r, \theta;\bar{r}, \bar{\theta})=\delta(\theta-\bar \theta) , & \theta, \bar{ \theta} \in (-\pi, \pi].
\end{cases}
\end{equation}
Since $(1-r^2)^2>0$, we can derive the Poisson kernel related to the Dirichlet problem (\ref{gfgf}) from the Euclidean case:
\begin{equation} \label{oio}
u( r,{\theta};\bar{r},\bar{\theta})=\frac{1}{2 \pi}  \frac{\bar{r}^2-r^2}{\bar{r}^2+r^2-2 r \bar{r}  \cos(\theta-\bar{\theta})}.
\end{equation}
Alternatively it is possible to obtain formula (\ref{oio}) from the Poisson kernel in $\mathbb{H}^2$ with a change of coordinates. In fact, in view of (\ref{ea}) and (\ref{rr}), formula (\ref{oio}) immediately follows. 

The Poisson kernel in (\ref{oio}) represents the law of the position occupied by the hyperbolic Brownian motion $\{B_{\mathbb{D}^2}(t): t \ge 0\}$ on $\mathbb{D}^2$ starting from $Q=(r,\theta) \in \mathbb{D}^2$ when it hits for the first time the boundary $\partial U$. We have 
$$\mathbb{P}_Q\{B_{\mathbb{D}^2}(T_{\bar{r}}) \in \mathrm{d} \bar{\theta}\}=\frac{1}{2 \pi}  \frac{\bar{r}^2-r^2}{\bar{r}^2+r^2-2 r \bar{r}  \cos(\theta-\bar{\theta})} \mathrm{d} \bar{\theta}.$$
We note that for $\bar{r} \to 1$ we have
\begin{equation} \label{lum}
\tilde{u}(r,\theta;\bar{\theta}):= \lim_{\bar{r} \to 1} u(r,{\theta};\bar{r},\bar{\theta}) =\frac{1}{2 \pi}  \frac{1-r^2}{1+r^2-2 r  \cos(\theta-\bar{\theta})}, \hspace{1cm} r<1.
 \end{equation}
Again  $\tilde{u}(r,\theta;\bar{\theta})$ represents the law of the position occupied by the hyperbolic Brownian motion in $\mathbb{D}^2$ when it hits for the first time the boundary of the hyperbolic disc $U$ with hyperbolic radius that goes to infinity. Result (\ref{lum}) is stated in Helgason \cite{helgason} page 34. For the $n$-dimensional case see Byczkowski and Malecki \cite{byczkowski2} formula (16). 

In view of (\ref{rr}), we can write (\ref{lum}) in hyperbolic coordinates as follows 
\begin{eqnarray*}
\frac{1}{2 \pi}\frac{1-r^2}{1+r^2-2 r \cos (\theta-\bar{\theta})}&=&\frac 1 {2 \pi}+ \frac 1 \pi  \sum_{n=1}^\infty  r^n  \cos n( \theta-\bar{\theta})=\frac 1 {2 \pi}+ \frac 1 \pi  \sum_{n=1}^\infty  \left(\frac{\cosh \eta-1}{\sinh \eta}   \right)^n  \cos n(\theta-\bar{\theta})\\
&=&\frac 1 {2 \pi} \frac{1}{\cosh \eta- \sinh \eta \cos (\theta-\bar{\theta})}, 
\end{eqnarray*}
which coincides with (\ref{lim}). On the other side it is well-known that under the conformal mapping (\ref{xxx}) the Poisson kernel (\ref{lum}) takes the form of the Cauchy distribution as it is shown in formula (\ref{fcv}). \\

Since $(1-r^2)^2>0$, the exit probabilities from the hyperbolic annulus $A=\{(r,\theta): r_1<r<r_2\}$ are easily derived from the Euclidean case. If the hyperbolic Brownian motion starts at $Q=(r,\theta) \in A$, we have 
\begin{equation} \label{wewe}
\mathbb{P}_Q\{T_{r_1}<T_{r_2}\}=\frac{\log r_2-\log r}{\log r_2-\log r_1}, \hspace{1cm} 0<r_1<r<r_2<1.
\end{equation}
Letting $r_2 \to 1$ in (\ref{wewe}) we obtain that
\begin{equation*}
\mathbb{P}_Q\{T_{r_1}<\infty\}=\frac{\log r}{\log r_1}<1, \hspace{1cm} 0<r_1<r<1.
\end{equation*}

\section{Brownian motion on the surface of a three-dimensional sphere}

The surface $S$ of  the unit-radius three dimensional sphere is a model of the elliptic geometry if geodesic lines are represented by great circles. We specify the position of an arbitrary point $p \in S$ with the couple $(\vartheta,\varphi)$  of spherical coordinates where $\vartheta \in [0,\pi]$ and $\varphi \in [0,2 \pi)$. 

If $U=\{(\vartheta,\varphi):   \vartheta>\bar{\vartheta}\}$ is the surface of a spherical cap on $S$ with center in the south pole, the Dirichlet problem on the surface of the sphere $S$ reads:
\begin{equation*} 
\begin{cases}
\left[ \frac{\partial^2}{\partial \vartheta^2} + \frac{1}{\tan \vartheta} \frac{\partial }{\partial \vartheta}+\frac{1}{\sin^2 \vartheta}  \frac{\partial^2}{\partial \varphi^2} \right] u( \vartheta, \varphi;\bar{\vartheta}, \bar{\varphi})=0, &0<\bar{\vartheta} < \vartheta < \pi, \\
u( \bar \vartheta, \varphi; \bar{\vartheta}, \bar{\varphi})=\delta({\varphi}-\bar{\varphi}) , & \varphi, \bar{ \varphi} \in [0, 2 \pi).
\end{cases}
\end{equation*}
Assuming that $u(\vartheta, \varphi; \bar{\vartheta}, \bar{\varphi})=T(\vartheta) F(\varphi)$ we immediately arrive at the following ordinary equations 
\begin{equation} \label{yah}
\begin{cases}
F''(\varphi)+\mu^2 \: F(\varphi)=0,\\
\sin^2 \vartheta \;T''(\vartheta)+\cos \vartheta \; \sin \vartheta\; T'(\vartheta)-\mu^2\; T(\vartheta)=0, 
\end{cases}
\end{equation} 
with $\mu \in \mathbb{R}$. With the change of variable $w=\cos \theta$, in the second equation of (\ref{yah}), we arrive at equation (\ref{ss}) with general solution (\ref{hhh}). Therefore, for $\mu=m \in \mathbb{N}$, the general solution to the second equation of (\ref{yah}) can be written as 
\begin{align} \label{xdx}
T(\vartheta)=C_1\left( \sqrt{\frac{1+\cos \vartheta}{1-\cos \vartheta}} \right)^m+C_2 \left( \sqrt{\frac{1-\cos \vartheta}{1+\cos \vartheta}} \right)^m.
\end{align}
We restrict ourselves to the increasing component of (\ref{xdx}) so that we have 
\begin{align*}
u( \vartheta, \varphi;\bar{\vartheta}, \bar{\varphi})=\sum_{m=0}^{\infty}  [A_m \cos(m \varphi)+B_m \sin(m \varphi)] \left( \sqrt{\frac{1-\cos \vartheta}{1+\cos \vartheta}} \right)^m.
\end{align*}
By imposing the boundary condition $u( \bar \vartheta, \varphi; \bar{\vartheta}, \bar{\varphi})=\delta({\varphi}-\bar{\varphi})$ and in view of (\ref{2.15}), we finally obtain that
\begin{align} \label{yuppy}
u(\vartheta, \varphi; \bar{\vartheta}, \bar{\varphi})&=\frac{1}{2 \pi}+\frac 1 \pi \sum_{m=1}^\infty  \cos(m (\varphi-\bar \varphi)) \left( \sqrt{\frac{1-\cos \vartheta}{1+\cos \vartheta}} \sqrt{  \frac{1+\cos \bar \vartheta}{1-\cos \bar  \vartheta}} \right)^m \nonumber \\
&=\frac{1}{2 \pi} \frac{1- \frac{1-\cos \vartheta}{1+\cos \vartheta}  \frac{1+\cos \bar \vartheta}{1-\cos \bar  \vartheta}}{1+\frac{1-\cos \vartheta}{1+\cos \vartheta}  \frac{1+\cos \bar \vartheta}{1-\cos \bar  \vartheta}-2 \sqrt{\frac{1-\cos \vartheta}{1+\cos \vartheta}  \frac{1+\cos \bar \vartheta}{1-\cos \bar  \vartheta}} \cos (\varphi- \bar \varphi)} \nonumber \\
&=\frac{1}{2 \pi} \frac{\cos   \vartheta- \cos  \bar \vartheta}{1- \cos \vartheta \cos \bar \vartheta- \sin  \vartheta \sin \bar \vartheta \cos (\varphi- \bar \varphi)}.
\end{align}

\begin{rem}
Since for spherical triangles the following Carnot formula holds
$$\cos \hat{\vartheta}=\cos \vartheta \cos \bar \vartheta+ \sin  \vartheta \sin \bar \vartheta \cos (\varphi- \bar \varphi),$$
we can rewrite (\ref{yuppy}) as follows 
\begin{equation} \label{lollo}u(\vartheta, \varphi; \bar{\vartheta}, \bar{\varphi})=\frac{1}{2 \pi} \frac{\cos \vartheta- \cos  \bar \vartheta}{1- \cos \hat{\vartheta}}.\end{equation}
\end{rem}

\begin{rem} We note that the Poisson kernel (\ref{yuppy}) is a proper probability law. In fact 
\begin{itemize}
\item In view of (\ref{lollo}) and observing that $ \bar \vartheta< \vartheta$ we have that (\ref{yuppy}) is positive.
\item Applying (\ref{omm}) with $a=1- \cos \vartheta \cos \bar  \vartheta$ and $b=- \sin \vartheta \sin \bar \vartheta$ we obtain that (\ref{yuppy}) integrates to one.  
\end{itemize}
For $\vartheta=0$ we obtain from  (\ref{yuppy}) the uniform law, while for $\bar \vartheta=\frac \pi 2$ we get $$u(\vartheta, \varphi; \frac \pi 2 , \bar{\varphi})=\frac{1}{2 \pi} \frac{\cos \vartheta}{1-\sin \vartheta \cos (\varphi- \bar \varphi)}.$$
\end{rem}

\begin{figure}[h]
\begin{minipage}[b]{7.85cm}
\centering
    \includegraphics[width=6.5cm, height=6.5cm]{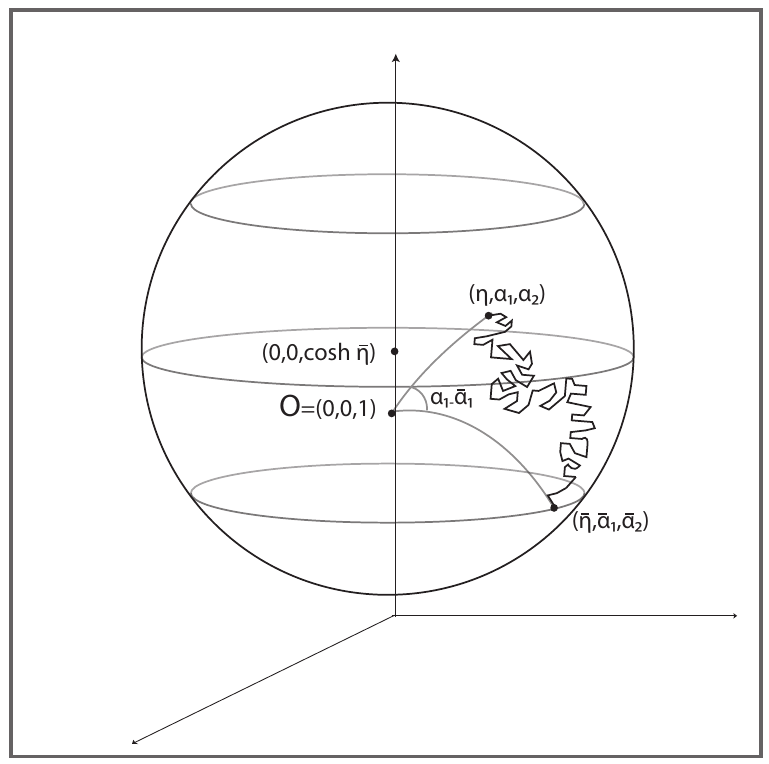}
    \caption{Brownian motion on $\mathbb{H}^3$ starting at $(\eta,\bm \alpha)$ and hitting the boundary of the hyperbolic ball} \label{sxx}
\end{minipage}
\hspace{0.25cm}
\begin{minipage}[b]{7.85cm}
\centering
    \includegraphics[width=6.5cm, height=6.5cm]{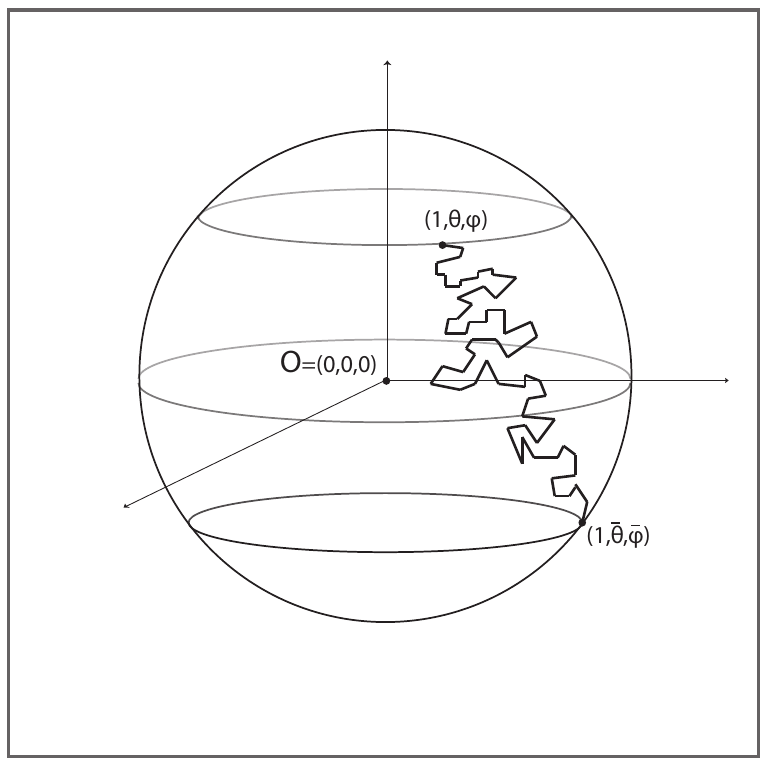}
    \caption{Spherical Brownian motion starting at $(\vartheta,\varphi)$ and hitting the boundary of the spherical disc} \label{sss}
\end{minipage}
 \end{figure}

\begin{rem}
Let $\{B_S(t): t \ge 0\}$ be a Brownian motion on the surface of the three dimensional sphere $S$ with starting point $p=(\vartheta,\varphi) \in S$ (see Figure \ref{sss}). The kernel in (\ref{yuppy}) represents the law of the position occupied by the spherical Brownian motion when it hits for the first time the boundary of the spherical cap $U$.
\end{rem}

In order to obtain the exit probabilities of $\{B_S(t): t \ge 0\}$ from a spherical annulus $A=\{(\vartheta,\varphi): \vartheta_2<\vartheta<\vartheta_1\}$ with center in the south pole of $S$, we consider the solution $v:(\vartheta_2,\vartheta_1) \to \mathbb{R}$ to the Laplace equation involving only the radial part 
$$\left[   \frac{\partial^2}{\partial \vartheta^2} + \frac{1}{\tan \vartheta} \frac{\partial }{\partial \vartheta}  \right] v(\vartheta)=0.$$
With the change of variable $w=\cos \vartheta$ we arrive at equation (\ref{xwxw}). With calculations analog to those performed in the proof of Theorem \ref{unabellalabel} we get 
\begin{equation} \label{roditore}
\mathbb{P}_p\{T_{\vartheta_1}<T_{\vartheta_2}\}=\frac{ \log \left| \frac{\cos \vartheta_2-1}{\cos \vartheta_2+1}  \right|- \log \left| \frac{\cos \vartheta-1}{\cos \vartheta+1}  \right|}{ \log \left| \frac{\cos \vartheta_2-1}{\cos \vartheta_2+1}  \right|- \log \left| \frac{\cos \vartheta_1-1}{\cos \vartheta_1+1}  \right|}, \hspace{1cm} \vartheta_2<\vartheta<\vartheta_1.
\end{equation} 
In particular for $\vartheta_2 \to \frac \pi 2$ formula (\ref{roditore}) reads
\begin{equation*}
\mathbb{P}_p\{T_{\vartheta_1}<T_{\frac \pi 2}\}=\frac{ \log \left| \frac{\cos \vartheta-1}{\cos \vartheta+1}  \right|}{ \log \left| \frac{\cos \vartheta_1-1}{\cos \vartheta_1+1}  \right|}=\frac{\log \left|\tan \frac \vartheta 2\right|}{\log \left|\tan \frac{ \vartheta_1} 2\right|}, \hspace{1cm} \frac \pi 2<\vartheta<\vartheta_1.
\end{equation*}

\begin{rem}
We note that replacing formally $\vartheta$ with $i \vartheta$ it is possible to extract from (\ref{yuppy}) and (\ref{roditore}) the Poisson kernel and the exit probabilities obtained for the hyperbolic plane, namely (\ref{statement}) and (\ref{uku}). This is because $\mathbb{H}^2$ can be viewed formally as a sphere with imaginary radius. For small values of $\theta$ we obtain instead results analogous to these obtained for the  Euclidean Brownian motion. In fact, for sufficiently small domains, Euclidean geometry is in force.  
 \end{rem}

\end{document}